\magnification=1200
\hoffset3cm

\input amstex

\catcode`\@=11
\font\myfont=cmss5
\def\myifdefined#1#2#3{%
  \expandafter\ifx\csname #1\endcsname \relax#2\else#3
    \fi}
\newcount\fornumber\newcount\artnumber\newcount\tnumber
\newcount\secnumber
\newcount\constnumber
\newcount\theonumber
\def\Ref#1{%
  \expandafter\ifx\csname mcw#1\endcsname \relax
    \warning{\string\Ref\string{\string#1\string}?}%
    \hbox{$???$}%
  \else \csname mcw#1\endcsname \fi}
\def\Refpage#1{%
  \expandafter\ifx\csname dw#1\endcsname \relax
    \warning{\string\Refpage\string{\string#1\string}?}%
    \hbox{$???$}%
  \else \csname dw#1\endcsname \fi}

\def\warning#1{\immediate\write16{%
            -- warning -- #1}}
\def\CrossWord#1#2#3{%
  \def\x{}%
  \def\y{#2}%
  \ifx \x\y \def\z{#3}\else
            \def\z{#2}\fi
  \expandafter\edef\csname mcw#1\endcsname{\z}
\expandafter\edef\csname dw#1\endcsname{#3}}
\def\Talg#1#2{\begingroup
  \edef\mmhead{\string\CrossWord{#1}{#2}}%
  \def\writeref{\write\refout}%
  \expandafter \expandafter \expandafter
  \writeref\expandafter{\mmhead{\the\pageno}}%
\endgroup}

\openin15=\jobname.ref
\ifeof15 \immediate\write16{No file \jobname.ref}%
\else      \input \jobname.ref \fi \closein15
\newwrite\refout
\openout\refout=\jobname.ref

\def\dfs#1{\myifdefined{mcwx#1}{}{\warning{multiply defined label #1}}
     \global\advance\secnumber by 1\global\fornumber=0\global\tnumber=0
      \the\secnumber.\if\the\theonumber0
      {\myfont #1}%
\else\relax\fi
     \Talg{#1}{\the\secnumber}
 \expandafter\gdef\csname mcwx#1\endcsname{a}\ignorespaces
}%

\def\dff#1{\myifdefined{mcwx#1}{}{\warning{multiply defined label #1}}
     \global\advance\fornumber by 1 \if\the\theonumber0 \text{\myfont #1}%
\else\relax\fi
\unskip (\the\secnumber.\the\fornumber)
\enspace
     \Talg{#1}{(\the\secnumber.\the\fornumber)}
  \expandafter\gdef\csname mcwx#1\endcsname{a} \ignorespaces
}%
\def\dfc#1{\myifdefined{mcwx#1}{}{\warning{multiply defined label #1}}
     \global\advance\constnumber by 1 \if\the\theonumber0 \text{\myfont #1}%
\else\relax\fi
\unskip \the\constnumber
     \Talg{#1}{\the\constnumber}
  \expandafter\gdef\csname mcwx#1\endcsname{a} \ignorespaces
}%


\def\dft#1{\myifdefined{mcwx#1}{}{\warning{multiply defined label #1}}%
     \global\advance\tnumber by 1%
       \if\the\theonumber0%
 \ignorespaces    {\myfont #1 }%
\else\relax\fi%
 \the\secnumber.\the\tnumber
    \Talg{#1}{\the\secnumber.\the\tnumber}
  \expandafter\gdef\csname mcwx#1\endcsname{a} \unskip\ignorespaces
}%
\def\dfro#1{\myifdefined{mcwx#1}{}{\warning{multiply defined label #1}}%
     \if\the\theonumber0%
      {\par\noindent\llap{\myfont #1}\par}%
\else\relax\fi%
     \Talg{#1}{(\the\rostercount@)}%
  \expandafter\gdef\csname mcwx#1\endcsname{a} \unskip\ignorespaces
}%
\def\dfa#1{\myifdefined{mcwx#1}{}{\warning{multiply defined label #1}}
     \global\advance\artnumber by 1
      \unskip      \the\artnumber%
     \Talg{art.#1}{\the\artnumber}
     \expandafter\gdef\csname mcwx#1\endcsname{a} \unskip\ignorespaces
 }%

\def\rf#1{\ifmmode\Ref{#1}\else $\Ref{#1}$\fi}

\def\donotshowtheoremlabels{\theonumber=1}

\def\rfa#1{\Ref{art.#1}}

\csname amsppt.sty\endcsname

\documentstyle{amsppt}

\vsize 18.8cm

\def\Li{\roman L}


\def\df#1{\leqno{\myifdefined{mcwx#1}{}{\warning{multiply
defined label #1}}\unskip
     \global\advance\fornumber by 1 \if\the\formunumber0 {\rm ref\{#1\}}%
\else\relax\fi
     \unskip (\the\secnumber.\the\fornumber)
    \Talg{#1}{(\the\secnumber.\the\fornumber)}\unskip
   \expandafter\gdef\csname mcwx#1\endcsname{a}\ignorespaces
}}%

\def\dfr#1#2{\myifdefined{mcwx#1}{}{\warning{multiply
defined label #1}}
     \Talg{#1}{#2}
  \expandafter\gdef\csname mcwx#1\endcsname{a} \unskip\ignorespaces
}%

\newdimen\mydim
\def\mycenter#1{\mydim=\hsize \advance \mydim by -50pt
\vtop{\hsize=\mydim
 \noindent\ignorespaces  #1 }}

\def\proof{\demo{Proof}} \def\eproof{\qed\enddemo}
\def\co{\colon}
\def\supp{\operatorname{\roman supp}}
\def\N{{\Bbb N}}
\def\Z{{\Bbb Z}}

\def\R{{\Bbb R}}
\def\D{\roman{d}}
\def\D{\roman{d}}
\def\eps{\varepsilon}
\def\is#1..#2..{\langle #1,#2\rangle}
\def\cici#1..#2..{\mathinner{[{#1}\ldotp\ldotp{#2}]}}
\def\ci#1..{\mathinner{[{#1}]}}
\def\oi#1..{\mathinner{]{#1}[}}
\def\ro#1..{\mathinner{[{#1}[}}
\def\lo#1..{\mathinner{]{#1}]}}
\def\ici#1..{\mathinner{[\![{#1}]\!]}} \def\Om{\Omega}
\donotshowtheoremlabels 
\nologo
\hoffset1.5cm

\topmatter
\def\mytime{\the\day.\the\month.\the\year-\the\time}
\title \ Attractors for damped  hyperbolic equations on
arbitrary unbounded domains
\endtitle

\author  Martino Prizzi --- Krzysztof P.
Rybakowski
\endauthor
\leftheadtext{\ M. Prizzi --- K. P.
Rybakowski}
\rightheadtext{\ Damped hyperbolic equations}

\address Martino Prizzi, Universit\`a degli Studi di
Trieste, Dipartimento di Matematica e Informatica, Via
Valerio, 12, 34127 Trieste, ITALY
\endaddress
\email
prizzi\@ dsm.univ.trieste.it
\endemail

\address Krzysztof P. Rybakowski, Universit\"at Rostock,
Institut f\"ur Mathematik, Universit\"atsplatz 1, 18055
Rostock, GERMANY
\endaddress
\email
krzysztof.rybakowski\@ mathematik.uni-rostock.de
\endemail
\date \enddate
\abstract We prove existence of global attractors for
damped hyperbolic equations of the form $$\aligned \eps
u_{tt}+\alpha(x) u_t+\beta(x)u- \sum_{ij}(a_{ij}(x)
u_{x_j})_{x_i}&=f(x,u),\quad x\in \Omega,\,t\in\ro0,\infty..,\\
u(x,t)&=0,\quad x\in \partial \Omega,\, t\in\ro0,\infty...
\endaligned$$ on an unbounded domain $\Omega$, without
smoothness assumptions on $\beta(\cdot)$,
$a_{ij}(\cdot)$, $f(\cdot,u)$ and
$\partial\Omega$, and $f(x,\cdot)$ having critical or subcritical growth.
\endabstract

\endtopmatter

\document \hoffset0cm

\head
\dfs{sec:intr}Introduction
\endhead

In this paper we study the existence of  global attractors
for  semilinear damped wave equations of the form
$$
 \aligned\eps u_{tt}+\alpha(x)
  u_t+\beta(x)u- \Li u&=f(x,u),\quad x\in
  \Omega,\,t\in\ro0,\infty..,\\ u(x,t)&=0,\quad x\in \partial
  \Omega,\, t\in\ro0,\infty...
 \endaligned\leqno\dff{261105-0951}
$$%
Here, $N\in \N$ and $\Omega$ is an {\it arbitrary\/} open
set in $\R^N$, bounded or not, $\eps>0$ is a constant
parameter, $\alpha$, $\beta\co\Omega\to\R$ and
$f\co\Omega\times\R\to \R$ are given functions and $\Li u
:=\sum_{ij}\partial_i(a_{ij}(x)\partial_j u)$ is a  linear
second-order differential operator in divergence form.

For  bounded domains $\Omega$  there are many results
concerning  the existence of attractors of~\rf{261105-0951}
under various assumptions on $\eps$, $\alpha$, $\beta$,
$\Li$ and $f$, including the pioneering works by Babin and
Vishik~\cite{\rfa{BV}}, Ghidaglia and Temam~\cite{\rfa{GT}}
and Hale and Raugel~\cite{\rfa{HR}}.

The unbounded domain case $\Omega=\R^3$ was considered
 in the important
papers~\cite{\rfa{F}, \rfa{F1}} by Feiresl.

In this paper  we assume that $\alpha\in L^\infty(\Omega)$,
$\alpha$ is bounded below by a positive constant and $\Li$
is  uniformly elliptic with coefficients functions lying in
$L^\infty(\Omega)$. We also assume that $\beta\in
L^p_{\roman u}(\R^N)$ with  $p>\max(1,N/2)$ and
$$
 \lambda_1=\inf\{\,E(u) \mid u\in
 H^1_0(\Omega),\,|u|^2_{L^2(\Omega)}=1\,\}>0\leqno\dff{171205-1355}
$$
where
$$
 E(u)=\int_{\Omega}\Bigl(\sum_{i,j=1}^N
 a_{ij}(x)\partial_iu(x)\partial_ju(x)+\beta(x)|u(x)|^2\Bigr)\,\D x.
$$
Here we denote by  $L^p_{\roman u}(\R^N)$  the set
of  measurable functions $v\co \R^N\to \R$ such that
$$
 |v|_{L^p_{\roman u}}:=\sup_{y\in \R^N}\left(\int_{
 B(y)}|v(x)|^p\,\D x\right)^{1/p}<\infty,
$$
where, for
$y\in\R^N$, $B(y)$ is the open unit cube in $\R^N$ centered
at $y$, cf ~\cite{\rfa{ACDR}}.

 We assume that the nonlinearity $f\co
\Omega\times \R\to\R$, $(x,u)\mapsto f(x,u)$ is measurable
in $x$, continuously differentiable in $u$ and satisfies
the growth assumptions $f(\cdot,0)\in L^2(\Omega) $ and
$$
 \text{$|\partial_u f(x,u)|\le \overline
 C(a(x)+|u|^{\overline\rho})$ for a.e. $x\in \Omega$ and
 every $u\in\R$.}
$$
Here $\overline C\ge 0$ and
$\overline\rho\ge0$ are constants with $2(\overline
\rho+1)\le 2^*:=(2N)/(N-2)$ for $N\ge 3$. If  $N\le 2$ or else if $N\ge 3$ and $2(\overline \rho+1)<2^*$, 
then $\overline\rho$ is called {\it subcritical\/}. If $N\ge 3$ and $2(\overline\rho+1)=2^*$, then
$\overline\rho$ is called {\it critical\/}.

In the subcritical case we also assume that  $a\in L^r_{\roman u}(\R^N)$ for some $r>\max(N,2)$, 
while in the critical case we assume that $a\in L^r(\Omega)+L^\infty(\Omega)$ for some $r\ge N$ and
$\alpha\in C^1(\Omega)$ with bounded derivatives. (Actually, our assumptions  concerning the functions
$\alpha$, $\beta$ and $a$ are somewhat more general than those listed above.)

Letting $F(x,u):=\int_0^uf(x,s)\,\D s $,
$(x,u)\in\Omega\times \R$, we  assume the dissipativity
conditions
$$
 \text{$f(x,u)u-\overline\mu F(x,u)\le c(x)$
 and $F(x,u)\le c(x)$ for a.e.  $x\in \Omega$ and every
 $u\in\R$}\leqno\dff{171205-1540}
$$ where $\overline\mu>0$
is a constant and  $c\in L^1(\Omega)$.

The goal of this paper is to prove that under the above hypotheses,  Equation~\rf{261105-0951} regarded 
as a system in $(u,v)$ where $v=u_t$, generates a nonlinear
continuous semigroup i.e. a {\it semiflow\/} $\pi_f$ on $Z=H^1_0(\Omega)\times L^2(\Omega)$ which has a 
global attractor.

Although our results hold for arbitrary open
sets $\Omega$, the emphasis here is on unbounded domains.

Condition~\rf{171205-1355} roughly means that the ground
state of the stationary Schr\"o\-dinger equation $$-\Li
u+\beta(x)u=0$$ on $\Omega$ with potential $\beta$ and with
Dirichlet boundary condition has positive energy.  In the
special case of $\beta\in L^1(\Omega)+ L^\infty(\Omega)$
with $\beta\ge 0$, condition~\rf{171205-1355} is equivalent
to the condition that $\int_G\beta(x)\,\D x=\infty$ for
every domain $G\subset\Omega$ that contains arbitrary large
balls. This was proved in~\cite{\rfa{AB1}, \rfa{AB2}}.

The dissipativity condition~\rf{171205-1540} was introduced
by Ghidaglia and Temam~\cite{\rfa{GT}} for the bounded
domain case. It is satisfied e.g. if there are constants
$\gamma$, $\nu\in\oi1,\infty..$ and a strictly positive
function $D\in L^1(\Omega)$ such that $F(x,u)\le D(x)$ for
all $x\in\Omega$, $u\in\R$  and the function $u\mapsto
(\gamma D(x)-F(x,u))^\nu$ is convex for a.e. $x\in\Omega$.

The proofs of our main results are based on Theorem~\rf{290505-2035} below, which provides
the so-called {\it tail estimates\/} for the solutions
$(u(t,x),u_t(t,x))$ of Equation~\rf{261105-0951}. For $\overline\rho$ subcritical, 
Theorem~\rf{290505-2035} implies that the semiflow $\pi_f$ is asymptotically compact on the phase space $Z$ 
(Lemma~\rf{050306-2033}) and this  proves the existence of a global attractor in the subcritical case
(Theorems~\rf{231105-1037}). For $\overline \rho$ critical we first use Theorem~\rf{290505-2035}
to show that $\pi_f$ is asymptotically compact with respect to the topology of the 
space $Y=L^2(\Omega)\times H^{-1}(\Omega)$  (Lemma~\rf{050306-2033}). Then we apply a method originally due
to J. Ball \cite{\rfa{Ba}} and elaborated by I. Moise, R. Rosa, X. Wang \cite{\rfa{MRW}} and G.
Raugel~\cite{\rfa{Ra}} to prove that $\pi_f$ is asymptotically compact on $Z$ (Theorem~\rf{070306-0658}).
This  establishes the existence of a global attractor in the critical case, see Theorem~\rf{070306-1631}.

The method of tail estimates was introduced by Wang~\cite{\rfa{W}} for
parabolic equations on unbounded domains and it was used by Fall and You~\cite{\rfa{FY}} to
establish the existence of an attractor of~\rf{261105-0951}
in the special case $\Omega=\R^N$, $\eps=1$,
$\beta(x)\equiv1$, $\Li=\Delta$, $\alpha(x)\equiv \lambda$
with $1\le \lambda <2$, and $f$  dissipative,  of {\it
sublinear\/} growth and having the special form
$f(x,u)=g(x)+\phi(u)$ with $g\in L^2(\R^N)$.

We should note that  our tail estimates  for the solution
component $u(t,x)$ do not depend in any way on the finite propagation speed property and are uniform in  the
parameter $\eps>0$. This allows us to prove singular semicontinuity results for
the family of attractors of Equation~\rf{261105-0951} as
$\eps\to 0$, cf. the forthcoming
publication~\cite{\rfa{PR}}.

For $N=3$ the exponent $\overline \rho$ is critical if $\overline\rho=2$ and subcritical 
if $\overline\rho<2$. In particular, Theorem~\rf{070306-1631} extends earlier results by
Feireisl~\cite{\rfa{F}}. 

In~\cite{\rfa{F1}} Feireisl proves existence of
attractors even in  the supercritical case
$2<\overline\rho<4$.
On the other hand, the
arguments in~\cite{\rfa{F}, \rfa{F1}}  require  additional
smoothness assumptions on $f(x,u)$ with respect to {\it
all\/} variables $(x,u)$ and some growth assumptions on
$|\partial_uf(x,0)|$ and $|\partial_xf(x,0)|$, while we do
not need any such condition here. Moreover, only the case
$\Omega=\R^3$ and $\Li=\Delta $ is considered
in~\cite{\rfa{F}, \rfa{F1}} and though the proofs do extend
to more general domains $\Omega$ and to more general
differential operators $\Li$, restrictions that  have to
imposed  are more stringent than the ones considered here.
In fact, the finite propagation speed property used in~\cite{\rfa{F}, \rfa{F1}} requires
some smoothness assumptions to be imposed on the
coefficient functions $a_{ij}(x)$ and on the boundary of
$\Omega$, cf.~\cite{\rfa{La1}}, while the Strichartz
estimates used in~\cite{\rfa{F1}} put some additional
restrictions both on the shape of $\Omega$ and  on the
coefficient functions $a_{ij}(x)$, cf.~\cite{\rfa{SS}}
and~\cite{\rfa{M}}.

This paper is organized as follows: in
Section~\rf{100405-1853} we collect some preliminary
concepts and results concerning semiflows, attractors and $(C_0)$-semigroups of linear operators.
We also establish an abstract differentiability result,
Theorem~\rf{p:diff}, which can frequently be used to
rigorously justify formal derivative calculations of
functionals along solutions of evolution equations. In
Section~\rf{150505-1251} we establish some general
estimates for linear damped wave equations and prove some
continuity and differentiability properties of Nemitski
operators. Finally, in Section~\rf{100405-1955}, we prove
our tail estimates and, as a consequence, establish the
existence of a global attractor of
Equation~\rf{261105-0951}.
\subhead Notation\endsubhead
For $a$ and $b\in\Z$ we write $\cici a..b..$ to denote the
set of all $m\in\Z$ with $a\le m\le b$.

Let $N\in\N$ be arbitrary.
Given a subset $S$ of $\R^N$
and a function $v\co S\to \R$ we denote by  $\tilde
v\co\R^N\to \R$  the trivial extension of $v$ defined by
$\tilde v(x)=0$ for $x\in \R^N\setminus S$.

Now let  $\Omega$ be
an arbitrary open set in $\R^N$. Given any measurable
function $v\co \Omega\to \R$ and any $p\in\ro1,\infty..$ we
set, as usual,
$$
 |v|_{L^p}=|v|_{L^p(\Omega)}:=\left(\int_\Omega|v(x)|^p\,\D
 x\right)^{1/p}\le\infty.
$$
Moreover, for $v\in
H^1_0(\Omega)$ we set
$|v|_{H^1}=|v|_{H^1(\Omega)}:=(|\nabla
u|_{L^2}^2+|u|_{L^2}^2)^{1/2}$.

If $k\in\N$ and $f$,  $g\co\Omega\to\R^k$ are such that
$\sum_{i=1}^kf_ig_i\in L^1(\Omega)$ then we set
$$\langle
f,g\rangle:=\int_{\Omega}\sum_{i=1}^kf_i(x)g_i(x)\,\D x.$$
We also use the common notation ${\Cal D}(\Omega)$ resp. ${\Cal D} '(\Omega)$ to denote the space of all test functions on $\Omega$, resp. all distributions on $\Omega$. If $w\in{\Cal D}'(\Omega)$ and $\varphi\in {\Cal D}(\Omega)$, then we use the usual functional notation $w(\varphi)$ to denote the value of $w$ at $\varphi$.

Given  a function $g\co
\Omega\times \R\to \R$, we denote by $\hat g$  the ({\it Nemitski\/}) operator which
associates with every  function $u\co \Omega\to \R$ the function
$\hat g(u)\co \Omega\to \R$ defined by
 $$
  \hat g(u)(x)= g(x,u(x)),\quad x\in \Omega.
 $$
All linear spaces considered in this paper are over the
real numbers.
\head\dfs{100405-1853}Preliminaries and an abstract
differentiability result
\endhead
We assume the reader's familiarity with attractor theory on metric spaces as expounded
in e.g.~\cite{\rfa{Ha}},~\cite{\rfa{La}} or, more recently,
in~\cite{\rfa{DCh}} and we just collect here a few relevant
concepts from that theory.
\definition{Definition}
Let $X$ be a metric space.
Recall that a {\sl local semiflow $\pi$ on $X$\/} is, by
definition, a continuous map from an open subset $D$ of
$\ro0,\infty..\times X$ to $X$ such that, for every $x\in
X$ there is an $\omega_x=\omega_{\pi,x}\in\lo0,\infty..$
with the property that $(t,x)\in D$ if and only if $t\in
\ro0,\omega_x..$,  and such that (writing $x\pi
t:=\pi(t,x)$ for $(t,x)\in D$) $x\pi 0=x$ for $x\in X$ and
whenever $(t,x)\in D$ and $(s,x\pi t)\in D$ then
$(t+s,x)\in D$ and $x\pi(t+s)=(x\pi t)\pi s$. Given an
interval $I$ in $\R$, a map $\sigma\co I\to X$ is called a
{\sl solution (of $\pi$)\/} if whenever $t\in I$ and
$s\in\ro0,\infty..$ are such that $t+s\in I$, then
$\sigma(t)\pi s$ is defined and $\sigma(t)\pi
s=\sigma(t+s)$. If $I=\R$, then $\sigma$ is called a {\sl
full solution (of $\pi$)\/}. A subset $S$ of $X$ is called
{\sl ($\pi$-)invariant\/} if for every $x\in S$ there is a
full solution $\sigma$ with $\sigma(\R)\subset S$ and
$\sigma(0)=x$.

Given a local semiflow $\pi$ on $X$ and a subset  $N$ of
$X$, we say that {\sl $\pi$ does not explode in $N$\/} if
whenever $x\in X$ and $x\pi \ro0,\omega_x..\subset N$, then
$\omega_x=\infty$. A {\sl global semiflow\/} is a local
semiflow with $\omega_x=\infty$ for all $x\in X$.

Now let $\pi$ be a global semiflow on $X$. A subset $A$ of
$X$ is called a {\sl global attractor (rel. to $\pi$)\/} if
$A$ is  compact, invariant and if for every bounded set $B$
in $X$ and every  open neighborhood $U$ of $A$ there is a
$t_{B,U}\in \ro0,\infty..$ such that $x\pi t\in U$ for all
$x\in B$ and all $t\in \ro t_{B,U},\infty..$. It easily
follows that a global attractor, if it exists, is uniquely
determined.

A subset $B$ of $X$ is called {\sl ($\pi$-)ultimately
bounded} if there is a $t_B\in \ro0,\infty..$ such the set
$\{\,x\pi t\mid x\in B,\,t\in \ro t_B,\infty..\, \}$ is
bounded.

$\pi$ is called {\sl asymptotically compact} if whenever
$B\subset X$ is  ultimately bounded, $(x_n)_n$ is a
sequence in $B$ and $(t_n)_n$ is a sequence in
$\ro0,\infty..$ with $t_n\to \infty$ as $n\to \infty$, then
the sequence $(x_n\pi t_n)_n$ has a convergent subsequence.
\enddefinition
The following result is well-known:
\proclaim{Proposition~\dft{211105-1325}} A global semiflow
$\pi$ on a metric space $X$ has a global attractor if and only if the
following conditions are satisfied:
\roster
\item $\pi$ is
asymptotically compact;
\item every bounded subset of $X$ is
ultimately bounded;
\item there is a bounded set $B_0$ in
$X$ with the property that for every $x\in X$ there is a
$t_x\in \ro0,\infty..$ such that $x\pi t_x\in
B_0$.
\endroster
\endproclaim
\proof This is
just~\cite{\rfa{DCh}, Corollary~1.1.4 and
Proposition~1.1.3}.
\eproof
We require a few results from the general theory of $(C_0)$-semigroups of linear operators.
\proclaim{Proposition~\dft{260206-1440}}
 Let $Z$ be a Banach space and $T(t)$, $t\in\ro0,\infty..$ be a $(C_0)$-semigroup of linear operators on $Z$ with generator $B\co D(B)\to Z$. Then, for every $z\in D(B)$ there is a unique function $u\co \ro 0,\infty..\to D(B)$ which is continuously differentiable into $Z$, $u(0)=z$ and
 $$
  u'(t)=Bu(t),\quad t\in \ro0,\infty.. .
 $$
 $u$ is given by $u(t)=T(t)z$ for all $t\in\ro0,\infty..$.
\endproclaim
\proof This follows from~\cite{\rfa{Go}, proof of Theorem II.1.2}
\eproof
\proclaim{Proposition~\dft{260206-1506}}
 Let $Z$ and $Y$ be  Banach spaces and $S_Z(t)$, $t\in \ro0,\infty..$ (resp. $S_Y(t)$, $t\in\ro0, \infty..$) be a $(C_0)$-semigroup of linear operators on $Z$ (resp. on $Y$)  with generator $C_Z\co D(C_Z)\to Z$ (resp. $C_Y\co D(C_Y) \to Y$). Let $\nu\co Z\to Y$ be a bounded linear map with $\nu(D(C_Z))\subset D(C_Y)$. If $\nu C_Z z=C_Y(\nu z)$ for all $z\in D(C_Z)$, then $\nu S_Z(t)z= S_Y(t)(\nu z)$ for all $z\in Z$ and all $t\in \ro 0,\infty..$.
\endproclaim
\proof
 An application of Proposition~\rf{260206-1440} shows that
 $\nu S_Z(t)z= S_Y(t)(\nu z)$ for all $z\in D(C_Z)$ and all $t\in \ro 0,\infty..$. The general case follows by density.
\eproof
\proclaim{Proposition~\dft{260206-1520}}
 Let $Z$ be a Banach space, $S_Z(t)$, $t\in \ro0,\infty..$  be a $(C_0)$-semigroup of linear operators on $Z$ with generator $C_Z\co D(C_Z)\to Z$ and $Q\co Z\to Z$ be linear and bounded. Then the operator $C_Z+Q\co D(C_Z)\to Z$ generates a $(C_0)$-semigroup $T_Z(t)$, $t\in \ro0,\infty..$ of linear operators on $Z$. Moreover,
 $$
  T_Z(t)z=S_Z(t)z+\int_0^t S_Z(t-s)QT_Z(s)z\,\D s \leqno\dff{260206-1527}
 $$
 for all $z\in Z$ and $t\in \ro0,\infty..$.
\endproclaim
\proof
 The first assertion follows from~\cite{\rfa{Go}, Theorem~I.6.4}. For $z\in D(C_Z)=D(C_Z+Q)$ and $t\in \ro0,\infty..$ formula~\rf{260206-1527} is proved using Proposition~\rf{260206-1440} and~\cite{\rfa{Go}, proof of Theorem~II.1.3 (ii)}. The general case follows by density.
\eproof

\proclaim{Proposition~\dft{230505-1255}}
 Let $Z$ be a Banach space and $T(t)$, $t\in\ro0,\infty..$
 be a $(C_0)$-semigroup of linear operators on $Z$ with
 infinitesimal generator $B\co D(B)\subset Z\to Z$.
 Suppose that $\Phi\co Z\to Z$ is a map which is Lipschitzian on
 bounded subsets of $Z$. Then, for each $\zeta\in Z$  there
 is a maximal $\omega_\zeta=\omega_{B,\Phi,\zeta}\in
 \lo0,\infty..$ and a uniquely determined continuous map
 $z=z_\zeta\co \ro0,\omega_\zeta..\to Z$ such that
 $$
  z(t)=T(t)\zeta+\int_0^t T(t-s)\Phi(z(s))\,\D s,\quad t\in
  \ro0,\omega_\zeta...
  \leqno\dff{230505-0723}
 $$
 Writing $\zeta \Pi t:=z_\zeta(t)$ for $t\in \ro0,\omega(\zeta)..$
 we obtain a local semiflow $\Pi=\Pi_{B,\Phi}$ on $Z$  which
 does not explode in bounded subsets of $Z$.
\endproclaim
\proof
 This follows from~\cite{\rfa{CH}, proofs of Theorem~4.3.4 and
 Proposition~4.3.7}.
\eproof
In the remaining part of this section we will establish a result which can be used to
rigorously justify formal differentiation of various
functionals along (mild) solutions of semilinear evolution
equations.
\proclaim{Theorem~\dft{p:diff}}
 Let $Z$ be a Banach space and $T(t)$, $t\in\ro0,\infty..$
 be a $(C_0)$-semigroup of linear operators on $Z$ with
 infinitesimal generator $B\co D(B)\subset Z\to Z$.
 Let $U$ be open in $Z$, $Y$ be a normed space and $V\co U\to Y$ be a function which, as a map from $Z$ to $Y$, is continuous at
 each point of $U$ and Fr\'echet differentiable at each
 point of $U\cap D(B)$. Moreover, let $W\co U\times Z\to Y$
 be a function which, as a map from $Z\times Z$ to $Y$, is
 continuous  and  such that $DV(z)(Bz+w)=W(z,w)$ for $z\in
 U\cap D(B)$ and $w\in Z$. Let $\tau\in \oi0,\infty..$ and
 $I:=\ci0,\tau..$. Let $\bar z\in U$, $g\co I\to Z$ be
 continuous and $z$ be a map from $I$ to $U$ such that
 $$
  z(t)=T(t)\bar z+\int_0^t T(t-s)g(s)\,\D s, \quad t\in I.
 $$
 Then the map $V\circ z\co I\to Y$ is differentiable
 and
 $$
  (V\circ z)'(t)=W(z(t),g(t)),\quad t\in I.
 $$
\endproclaim
\proof
 For
 $z\in D(B)$ set $|z|_{D(B)}:=|z|_Z+|Bz|_Z$. Since $B$ is closed, this defines a
 complete norm on $D(B)$.
 For $h\in \oi0,\infty..$ and $t\in I$
 set $M_h:=\sup_{t\in \ci0,h..}|T(t)|_{\Cal L(Z,Z)}$ and
 $$
  g_h(t):=(1/h)\int_0^h T(s)g(t)\,\D s.
 $$
 It is well-known that $g_h(t)\in D(B)$ and $Bg_h(t)=(1/h)(T(h)g(t)-g(t))$. Thus $g_h\co I\to D(B)$ and the estimate
 $$
  \multline|g_h(t)-g_h(t')|_{D(B)}=\bigl|(1/h)\int_0^h
  T(s)(g(t)-g(t'))\,\D s\bigr|_Z\\+
  \bigl|(1/h)(T(h)(g(t)-g(t'))-(g(t)-g(t')))\bigr|_Z\\ \le
  M_h|g(t)-g(t')|_Z+(1/h)(M_h+1)|g(t)-g(t')|_Z
  \endmultline
 $$
 shows that $g_h$ is continuous into $D(B)$. Moreover, we
 claim that $g_h(t)\to g(t)$ in $Z$ as $h\to 0^+$, uniformly
 on $I$. In fact, otherwise there is an
 $\eps\in\oi0,\infty..$ and sequences $(h_m)_{m\in\N}$ in
 $\oi0,\infty..$ and $(t_m)_{m\in \N}$ in $I$ such that
 $h_m\to 0$, $t_m\to t\in I$ and $|g_{h_m}(t_m)-g(t_m)|_Z\ge
 \eps$ for all $m\in \N$. But
 $$
  |g_{h_m}(t_m)-g(t_m)|_Z\le |g_{h_m}(t_m)-g(t)|_Z+ |g(t_m)-g(t)|_Z.
 $$
 Moreover,
 $$
  \multline|g_{h_m}(t_m)-g(t)|_Z=|(1/h_m)
  \int_0^{h_m}(T(s)g(t_m)-g(t))\,\D s|_Z\\ \le
  |(1/h_m)\int_0^{h_m}T(s)(g(t_m)-g(t))\,\D s|_Z+
  |(1/h_m)\int_0^{h_m}(T(s)g(t)-g(t))\,\D s|_Z.
  \endmultline
 $$
 W.l.o.g. $h_m\le 1$ for all $m\in\N$ so
 $$
  |(1/h_m)\int_0^{h_m}T(s)(g(t_m)-g(t))\,\D s|_Z\le
  M_1|g(t_m)-g(t)|_Z\to 0.
 $$
 Since $T(s)g(t)-g(t)\to 0$ in $Z$ as $s\to 0^+$, it follows that
 $$
  |(1/h_m)\int_0^{h_m}(T(s)g(t)-g(t))\,\D s|_Z\to 0.
 $$
 Putting things together we see that $|g_{h_m}(t_m)-g(t_m)|_Z\to 0$, a contradiction, proving our claim. Since $D(B)$ is dense in $Z$ there is a sequence $(\bar z_m)_{m\in\N}$ in $D(B)$ which converges to $\bar z$ in $Z$. Since $U$ is open in $Z$ we may assume that $\bar z_m\in U\cap D(B)$ for all $m\in \N$. Choose a sequence
 $(h_m)_{m\in\N}$ in $\oi0,\infty..$ converging to zero. For
 $m\in \N$ and $t\in I$ set
 $$
  z_m(t)=T(t)\bar z_m+\int_0^t T(t-s)g_{h_m}(s)\,\D s.
 $$
 It is well-known that $z_m(t)\in D(B)$. Moreover, the map $z_m\co I\to D(B)$ is continuous
 into $D(B)$ and differentiable into $Z$ with
 $z_m'(t)=Bz_m(t)+g_{h_m}(t)$ for $t\in I$. Furthermore, by
 what we have proved so far, $z_m(t)\to z(t)$ in $Z$ as
 $m\to \infty$, uniformly on $I$. It follows that $z_m(t)\in
 U\cap D(B)$ for some $m_0\in\N$ and all $m\ge m_0$ and
 $t\in I$. Moreover, by our hypotheses and by what we have
 proved so far, $(V\circ z_m)(t)\to (V\circ z)(t)$ and
 $(V\circ
 z_m)'(t)=DV(z_m(t))(Bz_m(t)+g_{h_m}(t))=W(z_m(t),g_{h_m}(t))\to
 W(z(t),g(t))$ in $Y$ uniformly on $I$. Thus $V\circ z$ is
 differentiable into $Y$  and $(V\circ z)'(t)=W(z(t),g(t))$
 for $t\in I$. The theorem is proved.
\eproof

\head
 \dfs{150505-1251}
 Damped hyperbolic equations
\endhead

For the rest of this paper, $N\in \N$ and $\Om$ is an
arbitrary open subset of $\R^N$, bounded or not.

Consider the following
\proclaim{Hypothesis~\dft{221105-1756}}\roster
\item"(1)"\dfr{171205-1239}{(1)} $a_0$,
$a_1\in\oi0,\infty..$ are constants and $a_{ij}\co
\Omega\to \R$, $i$, $j\in\cici1..N..$ are functions in
$L^\infty(\Omega)$ such that $a_{ij}=a_{ji}$, $i$,
$j\in\cici1..N..$, and for every $\xi\in\R^N$ and a.e.
$x\in\Omega$, $a_0|\xi|^2\le \sum_{i,j=1}^N
a_{ij}(x)\xi_i\xi_j\le a_1|\xi|^2 $.
$A(x):=(a_{ij}(x))_{i,j=1}^N$,
$x\in\Omega$.\item"(2)"\dfr{171205-1240}{(2)} $\beta\co
\Omega\to \R$ is a measurable function with the property
that \itemitem{$(i)$}\dfr{171205-1245}{(2i)} for every
$\overline \eps\in \oi0,\infty..$ there is a
$C_{\overline\eps}\in \ro0,\infty..$ with
$\bigl||\beta|^{1/2} u\bigr|_{L^2}^2\le \overline \eps
|u|_{H^1}^2+C_{\overline\eps}|u|_{L^2}^2 $ for all $u\in
H^1_0(\Omega)$;\itemitem{$(ii)$} $\lambda_1:=\inf\{\,\is
A\nabla u..\nabla u.. +\langle \beta  u,u\rangle \mid u\in
H^1_0(\Om),\,|u|_{L^2}=1\,\}>0$.\endroster\endproclaim
\remark{Remark} Note that, under Hypothesis~\rf{221105-1756}
item~\rf{171205-1239}, $\is A\nabla u..\nabla u..$ is defined and
under Hypothesis~\rf{221105-1756} item~\rf{171205-1245}, $\langle \beta
u,u\rangle$ is
defined.\endremark

The following lemma contains a condition ensuring that
$\beta$ satisfies
Hypothesis~\rf{221105-1756} item~$(2i)$.
\proclaim{Lemma~\dft{161105-1152}} Let $p\in\oi1,\infty..$
and $\beta\co\Omega\to\R$ be such that  $\tilde\beta\in
L^p_{\roman u}(\R^N)$. \roster \item If $p\ge N/2$, then
there is a $C\in\ro0,\infty..$ such that
$$\bigl||\beta|^{1/2} u\bigr|_{L^2}\le C |u|_{H^1} $$ for
all $u\in H^1_0(\Omega)$.\item If $p>N/2$, then for every
$\overline \eps\in \oi0,\infty..$ there is a
$C_{\overline\eps}\in \ro0,\infty..$ with
$$\bigl||\beta|^{1/2} u\bigr|_{L^2}^2\le \overline \eps
|u|_{H^1}^2+C_{\overline\eps}|u|_{L^2}^2 $$ for all $u\in
H^1_0(\Omega)$.\endroster\endproclaim \proof There is a
family $(y_j)_{j\in\N}$ of points in $\R^N$ such that
$\R^N=\bigcup_{j\in\N}\overline{B(y_j)}$ and the sets
$\overline{B(y_j)}$, $j\in\N$, are pairwise
non-overlapping. Write $B_j=B(y_j)$, $j\in\N$. Let
$p'=p/(p-1)$. Since $p\ge N/2 $ we have  $2p'\le 2^*$ for
$N\ge 3$.  Let $M\in\oi0,\infty..$ be a bound of the
imbedding $H^1(B)\to L^{2p'}(B)$ where $B=B(0)$. Then, by
translation, $M$ is also a bound of the imbedding
$H^1(B(y))\to L^{2p'}(B(y))$ for any $y\in\R^N$. Let $u\in
H^1_0(\Omega)$ be arbitrary. Then
$$\aligned\int_{\Omega}&|\beta(x)u^2(x)|\,\D
x=\int_{\R^N}|\tilde\beta(x)\tilde u^2(x)|\,\D x=
\sum_{j\in\N }\int_{B_j}|\tilde\beta(x)\tilde u^2(x)|\,\D x
\\&\le\sum_{j\in\N}
\left(\int_{B_j}|\tilde\beta(x)|^{p}\,\D x
\right)^{1/p}\left(\int_{B_j}|\tilde u(x)|^{2p'}\,\D x
\right)^{1/p'}\\&\le |\tilde \beta|_{L^p_{\roman
u}}\sum_{j\in\N} \left(\int_{B_j}|\tilde u(x)|^{2p'}\,\D x
\right)^{1/p'}\le |\tilde \beta|_{L^p_{\roman u}}
M^2\sum_{j\in\N}\bigl|\tilde
u|_{B_j}\bigr|_{H^1(B_j)}^2\\&= |\tilde \beta|_{L^p_{\roman
u}} M^2\sum_{j\in\N}\int_{B_j}(|\nabla\tilde
u(x)|^2+|\tilde u(x)|^2 )\,\D x=M^2|\tilde
\beta|_{L^p_{\roman u}}|\tilde u|_{H^1(\R^N)}^2\\&=
M^2|\tilde \beta|_{L^p_{\roman u}}|
u|_{H^1(\Omega)}^2.\endaligned$$ This proves the first part
of lemma. If $p>N/2$ we may choose $q$ such that $2p'<q$,
and $q<2^*$ for $N\ge 3$. We may then interpolate between
$2$ and $q$ and so, for every $\overline
\eps\in\oi0,\infty..$ there is a $C_{\overline
\eps}\in\ro0,\infty..$, independent of $u$ such that for
all $j\in \N$ $$\aligned\left(\int_{B_j}|\tilde
u(x)|^{2p'}\,\D x \right)^{1/2p'}&\le \overline \eps
\left(\int_{B_j}|\tilde u(x)|^{q}\,\D x \right)^{1/q}+
C_{\overline \eps}\left(\int_{B_j}|\tilde u(x)|^{2}\,\D x
\right)^{1/2}\\&\le \overline \eps M'\bigl|\tilde
u|_{B_j}\bigr|_{H^1(B_j)}+ C_{\overline\eps}\bigl|\tilde
u|_{B_j}\bigr|_{L^2(B_j)}.\endaligned$$
Here $M'\in\oi0,\infty..$ is a bound of the imbedding $H^1_0(B(y_j))\to L^q(B(y_j))$ for every $j\in \N$. Hence
$$\aligned\left(\int_{B_j}|\tilde u(x)|^{2p'}\,\D x
\right)^{1/p'}&\le 2(\overline \eps
M')^2\int_{B_j}(|\nabla\tilde u(x)|^2+|\tilde u(x)|^2 )\,\D
x+ 2C_{\overline\eps}^2\int_{B_j}|\tilde u(x)|^2\,\D
x.\endaligned$$ Thus, by the above computation,
$$\aligned\int_{\Omega}&|\beta(x)u^2(x)|\,\D x\le |\tilde
\beta|_{L^p_{\roman u}}\sum_{j\in\N}
\left(\int_{B_j}|\tilde u(x)|^{2p'}\,\D x
\right)^{1/p'}\\&\le |\tilde \beta|_{L^p_{\roman
u}}\sum_{j\in\N}\left(2(\overline \eps
M')^2\int_{B_j}(|\nabla\tilde u(x)|^2+|\tilde u(x)|^2 )\,\D
x+ 2C_{\overline\eps}^2\int_{B_j}|\tilde u(x)|^2\,\D
x\right)\\&=|\beta|_{L^p_{\roman u}}2(\overline \eps
M')^2|u|_{H^1}^2+|\beta|_{L^p_{\roman
u}}2C_{\overline\eps}^2|u|_{L^2}^2.\endaligned$$ Now an
obvious change of notation completes the proof of
the second part of the lemma. \eproof %
\remark{Remark~\dft{270206-0701}}
 Under Hypothesis~\rf{221105-1756} item~$(1)$ let the
 operator $\Li\co H^1_0(\Omega)\to{\Cal D}'(\Omega)$ be
 defined by
 $$
  \Li u=\sum_{i,j=1}^N\partial_i(a_{ij}\partial_ju),\quad u\in H^1_0(\Omega).
 $$
 The definition of distributional derivatives implies that
 $$
  (\Li u-\beta u)(v)=-\langle A\nabla u,\nabla v\rangle-\langle \beta u,v\rangle,\quad u\in H^1_0(\Omega),\,v\in \Cal D(\Omega).
 \leqno\dff{260206}
 $$
 It follows by density that
 $$
  \aligned\langle (\Li u-\beta u),v\rangle&=-\langle A\nabla u, \nabla v\rangle-\langle \beta u,v\rangle\\&\text{ for $u$, $v\in H^1_0(\Omega)$ with $\Li u-\beta u\in L^2(\Omega)$.}\endaligned
  \leqno\dff{270206-0706}
 $$
\endremark \proclaim{Lemma~\dft{111105-1838}} Assume
Hypothesis~\rf{221105-1756}. If $\kappa\in\ro0,\lambda_1..$
is arbitrary and if $\overline\eps$ and $\rho$ are chosen
such that
  $\overline \eps\in\oi0,a_0..$, $\rho\in\oi0,1..$
and $c:=\min
\bigl(\rho(a_0-\overline\eps),(1-\rho)(\lambda_1-\kappa)-\rho(\overline
\eps+C_{\overline\eps}+\kappa )\bigr)>0$ then $$c(|\nabla
u|_{L^2}^2+|u|_{L^2}^2)\le \is A\nabla u..\nabla u..
+\langle \beta  u,u\rangle-\kappa\is u..u..\le C(|\nabla
u|_{L^2}^2+|u|_{L^2}^2),\,\, u\in H^1_0(\Omega)$$ where
$C:=\max(a_1+\overline\eps,\overline\eps+C_{\overline\eps})$.
\endproclaim \proof This is just a simple
computation.\eproof
\proclaim{Lemma~\dft{270206-0836}}
 Assume Hypothesis~\rf{221105-1756} and let $\eps\in\oi0,\infty..$
 be arbitrary. For $u$, $v\in H^1_0(\Om)$ define
 $$\langle u,v\rangle_1=(1/\eps)\langle A\nabla u,\nabla v\rangle+
 (1/\eps)\langle\beta u,v\rangle .\leqno\dff{080505-1517}$$
 $\langle \cdot,\cdot\rangle_1$ is a scalar product on $H^1_0(\Omega)$ and the norm defined by this
 scalar product is equivalent to the usual norm on
 $H^1_0(\Om)$.

 For every $u\in H^1_0(\Omega)$ the distribution $-(1/\eps)\Li u+(1/\eps)\beta u\in {\Cal D}'(\Omega)$ can be uniquely extended to a continuous linear function $f_u$ from $H^1_0(\Omega)$ to $\R$.
 The operator
 $$
  \Lambda\co H^1_0(\Omega)\to H^{-1}(\Omega):=(H^1_0(\Omega))',\quad u\mapsto f_u
 $$
 is an isomorphism of normed spaces.
 The assignment
 $$
  (f,g)\in H^{-1}(\Omega)\times H^{-1}(\Omega)\mapsto
  \langle f, g\rangle_{-1}:=\langle \Lambda^{-1}(f),\Lambda^{-1}(g)\rangle_{1}
 $$
 defines a scalar product on $H^{-1}(\Omega)$. The norm defined by this scalar product is equivalent to the usual (operator) norm on $H^{-1}(\Omega)$.
\endproclaim
\proof
 This follows from Lemma~\rf{111105-1838} and the Lax-Milgram theorem.
\eproof
\proclaim{Proposition~\dft{070505-1208}} Assume
Hypothesis~\rf{221105-1756} and let $\alpha_0$,
$\alpha_1\in\ro0,\infty..$ and $\eps\in\oi0,\infty..$ be
arbitrary. Let $\alpha\co\Omega\to\R $ be a measurable
function with $\alpha_0\le \alpha(x)\le \alpha_1$ for a.e.
$x\in\Omega$.
 Set $Z =H^1_0(\Om)\times
L^2(\Om)$ and endow $Z$ with the usual norm $|z|_Z$ defined
by $$|z|_{Z}^2=|\nabla
z_1|^2_{L^2}+|z_1|^2_{L^2}+|z_2|^2_{L^2},\quad
z=(z_1,z_2).$$ Define $D(B)=D(B_{\alpha,\beta,\eps})$ to be
the set of all $(z_1,z_2)\in Z$ such that $z_2\in
H^1_0(\Om)$ and   $\Li z_1-\beta z_1$ (in the distributional sense)
lies in $L^2(\Om)$. Let $B=B_{\alpha,\beta,\eps}\co D(B)\to
Z$ be defined by $$B(z_1,z_2)=(z_2,-(1/\eps)\alpha
z_2-(1/\eps)\beta z_1+(1/\eps)\Li z_1),\quad z=(z_1,z_2)\in
D(B).$$ Under these hypotheses, $B$ is the generator of a
$(C_0)$-semigroup $T(t)=T_{\alpha,\beta,\eps}(t)$,
$t\in\ro0,\infty..$ on $Z$.  If, in addition, $\alpha_0>0$,
then there are  real constants
$M=M(\alpha_0,\alpha_1,\eps,\lambda_1)>0$,
$\mu=\mu(\alpha_0,\alpha_1,\eps,\lambda_1)>0$ such that
$$|T(t)z|_Z\le M e^{-\mu t}|z|_Z,\quad z\in Z,\, t\in \ro
0,\infty...\leqno\dff{080505-1704}$$
\endproclaim
\proof
 On $Z$ define the scalar product
$$\langle\!\langle(u_1,u_2),(w_1,w_2)\rangle\!\rangle=\langle
u_1,w_1\rangle_1+\langle u_2,w_2\rangle_{L^2}.
\leqno\dff{080505-1655}$$ It follows from Lemma~\rf{270206-0836} that the norm
$\|(u_1,u_2)\|=
\langle\!\langle(u_1,u_2),(u_1,u_2)\rangle\!\rangle^{1/2}$
is equivalent to the norm $|(u_1,u_2)|_Z$.

Now, for $(z_1,z_2)\in D(B)$, we obtain using~\rf{270206-0706} $$\aligned
\langle\!\langle B(z_1,z_2),&(z_1,z_2)\rangle\!\rangle=\\&
\langle z_2,z_1 \rangle_1+\langle -(1/\eps)\alpha
z_2-(1/\eps)\beta z_1+(1/\eps) \Li z_1,z_2
\rangle=-(1/\eps)\is \alpha z_2..z_2... \endaligned$$ Thus
$B$ is dissipative by~\cite{\rfa{CH}, Proposition~2.4.2}.
Let us now show that $B$ is $m$-dissipative. We
use~\cite{\rfa{CH}, Proposition~2.2.6} and so we only need
to show that for every $\lambda\in\oi0,\infty..$ and for
every $(f,g)\in Z$ there is a $(z_1,z_2)\in D(B)$ with
$$(z_1,z_2)-\lambda B(z_1,z_2)=(f,g)
\leqno\dff{080505-1549}$$ Now~\rf{080505-1549} is
equivalent to the validity of the two equations
$$z_2=(1/\lambda)(z_1-f)\leqno\dff{080505-1555}$$ and
$$((1/\lambda)+(1/\eps)\alpha+(1/\eps)\lambda\beta)z_1-(1/\eps)\lambda\Li
z_1=
g+((1/\lambda)+(1/\eps)\alpha)f.\leqno\dff{080505-1559}$$
Lemma~\rf{111105-1838} and the Lax-Milgram theorem (cf
~\cite{\rfa{CH}, proof of Proposition~2.6.1}) imply that
equation~\rf{080505-1559} can be solved for $z_1\in
H^1_0(\Om)$ with $\Li z_1-\beta z_1\in L^{2}(\Om)$. Now
equation~\rf{080505-1555} can be solved for $z_2\in
H^1_0(\Om)$. It follows that, indeed, $B$ is
$m$-dissipative and so, by the Hille-Yosida-Phillips
theorem, $B$ generates a $(C_0)$-semigroup $T(t)$,
$t\in\ro0,\infty..$, of linear operators on $Z$.

Now suppose $\alpha_0>0$. Choose $\mu$ such that
$$0<2\mu\le\min( 1,\alpha_0/(2\eps),
\lambda_1/(\eps+\alpha_1)).\leqno\dff{150505-1308}$$ We now
prove that for every $(u_1,u_2)\in Z$
$$\|T(t)(u_1,u_2)\|\le 2e^{-\mu t}\|(u_1,u_2)\|,\quad
t\in\ro0,\infty...\leqno \dff{080505-1703}$$ This
proves~\rf{080505-1704} in view of the equivalences of the
two above norms on $Z$. By density, it is sufficient to
prove~\rf{080505-1703} for $(u_1,u_2)\in D(B)$. Therefore,
let $(u_1,u_2)\in D(B)$ be arbitrary and define
$(z_1(t),z_2(t))=T(t)(u_1,u_2)$, $t\in\ro0,\infty..$. Then
the map $t\mapsto z(t)=(z_1(t),z_2(t))$  is differentiable
into $Z$, $z(t)\in D(B)$ and $\dot z(t)=Bz(t)$ for $t\in
\ro0,\infty..$. For $t\in\ro0,\infty..$ let $$\aligned
w(t)&=4\mu\langle z_1(t),z_2(t) \rangle+\is
z_2(t)..z_2(t)..\\&+ 2(1/\eps)\mu\is \alpha
z_1(t)..z_1(t).. +(1/\eps)\langle\beta
z_1(t),z_1(t)\rangle+ (1/\eps)\is A\nabla z_1(t)..\nabla
z_1(t).. .\endaligned\leqno\dff{080505-1724}$$ It follows
that $w$ is differentiable and a simple calculation
 shows $$\aligned(1/2)\dot
w(t)&=\is(2\mu-(1/\eps)\alpha) z_2(t)..z_2(t)..\\&-2\mu
(1/\eps)\langle \beta z_1(t),z_1(t)\rangle-2 (1/\eps)\mu\is
A\nabla z_1(t)..\nabla z_1(t)..\\&
\le\is(2\mu-(1/\eps)\alpha_0) z_2(t)..z_2(t)..\\&-2\mu
(1/\eps)\langle \beta z_1(t),z_1(t)\rangle-2 (1/\eps)\mu\is
A\nabla z_1(t)..\nabla z_1(t)..
\endaligned\leqno\dff{080505-1732}$$ By~\rf{080505-1724}
$$\aligned w(t)&\le4\mu((1/2) \is z_1(t)..z_1(t)..
+(1/2)\is z_2(t)..z_2(t)..)+\is z_2(t)..z_2(t)..\\&+
2(1/\eps)\mu\is\alpha z_1(t)..z_1(t)..\\&+(1/\eps)\langle
\beta z_1(t),z_1(t)\rangle+ (1/\eps)\is A\nabla
z_1(t)..\nabla z_1(t)..\\&\le (2\mu+1)\is z_2(t)..z_2(t)..+
2\mu(1+(1/\eps)\alpha_1)\is
z_1(t)..z_1(t)..\\&+(1/\eps)\langle \beta
z_1(t),z_1(t)\rangle+ (1/\eps)\is A\nabla z_1(t)..\nabla
z_1(t)...\endaligned$$Now $$\aligned 2\|z(t)\|^2&=2\is
z_2(t)..z_2(t)..+ (1/\eps)\langle \beta
z_1(t),z_1(t)\rangle+ (1/\eps)\is A\nabla z_1(t)..\nabla
z_1(t)..\\&+(1/\eps)\langle \beta z_1(t),z_1(t)\rangle+
(1/\eps)\is A\nabla z_1(t)..\nabla z_1(t)..
\endaligned$$By~\rf{150505-1308} $$2\is z_2(t)..z_2(t)..\ge
(2\mu+1) \is z_2(t)..z_2(t)..$$and
$$\aligned(1/\eps)\langle \beta z_1(t),z_1(t)\rangle&+
(1/\eps)\is A\nabla z_1(t)..\nabla z_1(t)..\\&\ge(1/\eps)
\lambda_1\is z_1(t)..z_1(t)..\ge
2\mu(1+(1/\eps)\alpha_1)\is z_1(t)..z_1(t)...\endaligned
$$Putting things together we see that $$w(t)\le
2\|z(t)\|^2,\quad t\in
\ro0,\infty...\leqno\dff{080505-2209}$$ Moreover,
by~\rf{150505-1308} $$\aligned w(t)&\ge-4\mu((1/2)4\mu \is
z_1(t)..z_1(t)..+(1/2)(1/4\mu)\is z_2(t)..z_2(t)..) +\is
z_2(t)..z_2(t)..\\&+ 2(1/\eps)\mu\is\alpha
z_1(t)..z_1(t)..+(1/\eps)\langle \beta
z_1(t),z_1(t)\rangle+ (1/\eps)\is A\nabla z_1(t)..\nabla
z_1(t)..\\&\ge (1/2)\is z_2(t)..z_2(t)..+
2\mu((1/\eps)\alpha_0-4\mu)\is
z_1(t)..z_1(t)..\\&+(1/\eps)\langle \beta
z_1(t),z_1(t)\rangle+ (1/\eps)\is A\nabla z_1(t)..\nabla
z_1(t)..\ge(1/2)\|z(t)\|^2.\endaligned$$ Thus
 $$w(t)\ge(1/2)\|z(t)\|^2,\quad t\in\ro0,\infty...\leqno\dff{080505-2224}$$
By~\rf{080505-2209}, \rf{080505-1732} and~\rf{150505-1308}
$$\aligned\mu w(t)&\le 2\mu\|z(t)\|^2=2\mu\is
z_2(t)..z_2(t)..+ 2\mu(1/\eps)\langle \beta
z_1(t),z_1(t)\rangle\\&+ 2(1/\eps)\mu\is A\nabla
z_1(t)..\nabla z_1(t)..\le ((1/\eps)\alpha_0-2\mu)\is
z_2(t)..z_2(t)..\\&+ 2\mu(1/\eps)\langle \beta
z_1(t),z_1(t)\rangle+ 2(1/\eps)\mu\is A\nabla
z_1(t)..\nabla z_1(t)..\le-(1/2)\dot w(t) \endaligned$$ so
$$\dot w(t)\le -2\mu w(t),\quad t\in
\ro0,\infty...\leqno\dff{080505-2239}$$ \rf{080505-2209},
\rf{080505-2224} and~\rf{080505-2239} imply that
$$\|z(t)\|^2\le 4e^{-2\mu t}\|z(0)\|^2,\quad t\in
\ro0,\infty..$$ and this in turn implies~\rf{080505-1703}.
The theorem is proved. \eproof
\proclaim{Proposition~\dft{270206-1500}} Assume
 Hypothesis~\rf{221105-1756} and let $\eps\in\oi0,\infty..$ be arbitrary. Define $C_Z:=B_{\alpha,\beta,\eps}$ and   $S_Z(t):=T_{\alpha,\beta,\eps}(t)$, $t\in \ro0,\infty..$ with $\alpha\equiv0$.
 Moreover, let $Y=L^2(\Omega)\times H^{-1}(\Omega)$ and define the  operator $C_Y\co D(C_Y)\to Y$ by $D(C_Y)=H^1_0(\Omega)\times L^2(\Omega)$ and
 $$
  C_Y(z_1,z_2)=(z_2,-\Lambda(z_1))
 $$
 where $\Lambda$ is defined in Lemma~\rf{270206-0836}.
 $C_Y$ is the generator of a $(C_0)$-semigroup $S_Y(t)$, $t\in\ro0,\infty..$ of linear operators on $Y$.

 Finally,
 $$
  \nu S_Z(t)z=S_Y(t)(\nu z),\quad z\in Z,\,t\in\ro0,\infty..
 $$
 where $\nu\co Z\to Y$ is the inclusion map.
\endproclaim
\proof
 On $Y$ define the scalar product
 $$
  \langle\!\langle(u_1,u_2),(w_1,w_2)\rangle\!\rangle_Y=\langle
  u_1,w_1\rangle_{L^2}+\langle u_2,w_2\rangle_{-1}.
 $$
 It follows from Lemma~\rf{270206-0836} that the norm
 defined by this scalar product is equivalent to the usual norm on $Y$.
 Now, for $(y_1,y_2)\in D(C_Y)$, we easily obtain
 $$
  \langle\!\langle C_Y(y_1,y_2),(y_1,y_2)\rangle\!\rangle_Y=0.
 $$
 Thus
 $B_Y$ is dissipative.
 Using the same arguments as in the proof of Proposition~\rf{070505-1208} (with $\alpha\equiv0$)
 we can show that for every $\lambda\in\oi0,\infty..$ and for
 every $(f,g)\in Y$ there is a $(y_1,y_2)\in D(C_Y)$ with
 $
  (y_1,y_2)-\lambda C_Y(y_1,y_2)=(f,g)
 $. Thus $C_Y$ is $m$-dissipative and this proves the first assertion. Since, by the definitions of $C_Z$ and $C_Y$, $\nu D(C_Z)\subset D(C_Y)$ and $\nu C_Z(z_1,z_2)=C_Y\nu(z_1,z_2)$ for all $(z_1,z_2)\in D(C_Z)$, the second assertion follows from Proposition~\rf{260206-1506}.
\eproof
\proclaim{Proposition~\dft{040306-2207}}
 Let $\alpha$, $Z$  and $T(t)$ be as in Proposition~\rf{070505-1208} and $Y$ be as in Proposition~\rf{270206-1500}. Suppose that
 $$(\exists\,C_{\dfc{050306-1327}}\in\ro0,\infty..)(\forall\,z\in L^2(\Omega))\,|\alpha z|_{H^{-1}}\le C_{\rf{050306-1327}} |z|_{H^{-1}}.\leqno\dff{050306-0725}$$
 Then there are constants $C_{\dfc{050306-1135}}$ and $C_{\dfc{050306-1136}}\in\ro0,\infty..$ such that
 $$
  |T(t)z|_Y\le C_{\rf{050306-1135}}e^{C_{\rf{050306-1136}}t}|z|_Y,\quad t\in\ro0,\infty..,\,z\in Z.
 $$
\endproclaim
\proof
 Define the bounded  linear map $Q\co Z\to Z$ by $(z_1,z_2)\mapsto (0,-\alpha z_2)$.   By Propositions~\rf{260206-1520} and~\rf{260206-1506} we have, for $z\in Z$ and $t\in\ro0,\infty..$
 $$
     T(t)z=S_Z(t)z+\int_0^t S_Z(t-s)QT(s)z\,\D s
        =S_Y(t)z+\int_0^t S_Y(t-s)QT(s)z\,\D s.
 $$
 There are constants $C_{\dfc{050306-1329}}$ and $C_{\dfc{050306-1330}}\in\ro0,\infty..$ such that
 $$
  |S_Y(t)y|_Y\le C_{\rf{050306-1329}}e^{C_{\rf{050306-1330}}t}|y|_Y,\quad t\in\ro0,\infty..,\,y\in Y.
 $$
 Using~\rf{050306-0725} we now obtain, for $z\in Z$ and $t\in\ro0,\infty..$
 $$
  \aligned
  |T(t)z|_Y&\le|S_Y (t)z|_Y+\int_0^t |S_Y(t-s)QT(s)z|_Y\,\D s\\
           &\le C_{\rf{050306-1329}}e^{C_{\rf{050306-1330}}t}|z|_Y
           +\int_0^t C_{\rf{050306-1329}} e^{C_{\rf{050306-1330}}(t-s)}C_{\rf{050306-1327}}|T(s)z|_Y\,\D s.
  \endaligned
 $$
  Now Gronwall's lemma completes the proof.
\eproof
The next result provides a sufficient condition for the validity of~\rf{050306-0725}.
\proclaim{Lemma~\dft{160505-0804}} If $a\in C^1(\Om)\cap
W^{1,\infty}(\Om)$ and $u\in H^1_0(\Om)$, then $a u\in
H^1_0(\Om)$ and $\partial_i(a u)=(\partial_i a) u+
a\partial_i u$, $i\in \cici1..N..$. Moreover, $|a
u|_{H^1_0}\le (2N+1)^{1/2}|a|_{W^{1,\infty}}|u|_{H^1_0}$.
Furthermore,
$$
 |az|_{H^{-1}}\le (2N+1)^{1/2}|a|_{W^{1,\infty}}|z|_{H^{-1}},\quad z\in L^2(\Omega).$$
Finally, if $U$ is an open subset of $\Omega$ and $a|_U\in
C^1_0(U)$ then $(a u)|_U\in H^1_0(U)$. \endproclaim \proof
Set $u_{(i)}=(\partial_i a) u+a\partial_i u$, $i\in
\cici1..N..$. There is a sequence $(v_n)_{n\in\N}$ in
$C^1_0(\Om)$ converging to $u$ in $H^1(\Om)$. It follows
that $a v_n\in C^1_0(\Om)$ and $\partial_i(a
v_n)=(\partial_i a) v_n+a\partial_i v_n$ for $n\in \N$ and
$i\in \cici1..N..$. H\"older's inequality implies that, for
$\varphi\in C^1_0(\Om)$ and $i\in \cici1..N..$, $a v_n\to a
u$ and $\partial_i(a v_n)\to u_{(i)}$ in $L^2(\Om)$ while
$\varphi\partial_i(a v_n)\to \varphi u_{(i)}$ and
$(\partial_i\varphi) a v_n\to (\partial_i\varphi) a u$ in
$L^1(\Om)$. Since $\langle \varphi,\partial_i(a
v_n)\rangle_{L^2}=-\langle \partial_i\varphi,a
v_n\rangle_{L^2}$ for $n\in\N$ and $i\in\cici1..N..$, it
follows that $au\in H^1(\Om)$, $\partial_i (au)=u_{(i)}$
for all $i\in \cici1..N..$ and $$\lim_{n\to\infty}|(a
v_n)-(a u)|_{H^1}=0.\leqno\dff{160505-1022}$$ Thus $a u\in
H^1_0(\Om)$. This proves the first part of the lemma. It
follows that $$\aligned|a& u|_{H^1_0}^2=|a
u|_{L^2}^2+\sum_{i=1}^N |\partial_i(a u)|_{L^2}^2= |a
u|_{L^2}^2+\sum_{i=1}^N |(\partial_i a) u+a\partial_i
u|_{L^2}^2\\&\le
|a|_{W^{1,\infty}}^2|u|_{L^2}^2+\sum_{i=1}^N|a|_{W^{1,\infty}}^2
(|u|_{L^2}+|\partial_i u|_{L^2})^2\\&\le
|a|_{W^{1,\infty}}^2(|u|_{L^2}^2+\sum_{i=1}^N
(2|u|_{L^2}^2+2|\partial_i u|_{L^2}^2))=
|a|_{W^{1,\infty}}^2((2N-1)|u|_{L^2}^2+2| u|_{H^1_0}^2)
\\&\le |a|_{W^{1,\infty}}^2(2N+1)|
u|_{H^1_0}^2.\endaligned$$
If $z\in L^2(\Omega)$ then $az\in L^2(\Omega)$ and for $v\in H^1_0(\Omega)$ with $|v|_{H^1}\le 1$  we have $av\in H^1_0(\Omega)$ and
$$
 |\langle az,v\rangle|=|\langle z,av\rangle|\le |z|_{H^{-1}}|av|_{H^1}\le (2N+1)^{1/2}|a|_{W^{1,\infty}}|z|_{H^{-1}}.
$$ This proves the second and third part of
the lemma. Finally, if $a|_U\in C^1_0(U)$ then $(a
v_n)|_U\in C^1_0(U)$ for all $n\in \N$ and since,
by~\rf{160505-1022}, $(a v_n)|_U\to (a u)|_U$ in $H^1(U)$,
it follows that $(a u)|_U\in H^1_0(U)$. The lemma is
proved
\eproof
\proclaim{Proposition~\dft{150505-1711}} Let
$a\in C^1_0(\R^N)$ and $r\in \ro2,\infty..$ be arbitrary.
If $N\ge 3$, then assume also that $r< 2^*$. Under these
assumptions the map $h\co H^1_0(\Om)\to L^r(\Om)$,
$u\mapsto a|_{\Om}\cdot u$, is defined and is linear and
compact.\endproclaim\proof  There is an open ball $U$ in
$\R^N$ such that $\supp a\subset U$. Define the following
maps: $$\aligned &h_1\co H^1_0(\Om)\to
H^1_0(\R^N),\,u\mapsto \tilde u,\,\, h_2\co H^1_0(\R^N)\to
H^1_0(U),\,v\mapsto (a v)|_U,\\& h_3\co H^1_0(U)\to
L^r(U),\,v\mapsto v,\,\,h_4\co L^r(U)\to
L^r(\R^N),\,v\mapsto \tilde v \\& h_5\co L^r(\R^N)\to
L^r(\Om),\,v\mapsto v|_\Om.\endaligned $$ Clearly, the maps
$h_1$, $h_4$ and $h_5$ are  defined, linear and bounded,
the map $h_2$ is defined, linear and bounded in view of
Lemma~\rf{160505-0804} with $\Omega:=\R^N$, while $h_3$ is
defined, linear and compact by  Rellich embedding theorem.
Since, for all $u\in H^1_0(\Om)$, $(h_5\circ h_4\circ
h_3\circ h_2\circ h_1)(u)=a|_\Om\cdot u$, it follows that
$h$ is defined and $h=h_5\circ h_4\circ h_3\circ h_2\circ
h_1$ so $h$ is linear and compact. \eproof
 \definition{Definition} A
function $f\co \Om\times\R\to \R$, $(x,u)\mapsto f(x,u)$ is
said to satisfy a {\sl $C^0$- (resp. $C^1$-)Carath\'eodory
condition\/}, if for every $u\in \R$ the partial map
$x\mapsto f(x,u)$ is Lebesgue-measurable and for a.e. $x\in
\Om$ the partial map $u\mapsto f(x,u)$ is continuous (resp.
continuously differentiable).

If $f\co \Om\times\R\to \R$, $(x,u)\mapsto f(x,u)$
satisfies a  $C^0$-Carath\'eodory condition, define the
function  $F\co \Om\times \R\to \R$ by $$F(x,u)=\int_0^u
f(x,s)\,\D s,
 $$ whenever $s\mapsto f(x,s)$ is continuous and $F(x,u)=0$
 otherwise. $F$ is called the {\it canonical primitive of
 $f$\/}.

Given $\overline C$, $\overline\rho\in \ro0,\infty.. $, a
measurable function $a\co\Omega\to\R$ and a null set
$M\subset \Om$, a function $g\co (\Om\setminus
M)\times\R\to \R$, $(x,u)\mapsto g(x,u)$ is said to satisfy
a {\sl $(\overline C,\overline \rho,a)$-growth
condition\/}, if $| g(x,u)|\le \overline
C(|a(x)|+|u|^{\overline\rho})$ for every $x\in \Om\setminus
M$ and every $u\in\R$.
The number $\overline \rho$ is called {\it subcritical\/} if $N\le 2$ or ($N\ge 3$ and $\overline \rho<(2^*/2)-1$). $\overline \rho$ is called {\it critical\/} if $N\ge 3$ and $\overline \rho=(2^*/2)-1$.\enddefinition

\proclaim{Proposition~\dft{100405-0808}} Let $f$ satisfy a
$C^1$-Carath\'eodory condition and $\partial_u f$ satisfy a
$(\overline C,\overline \rho,a)$-growth condition. Let $F$
be the canonical primitive of $f$. Then, for a.e.
$x\in\Omega$ and all $u$, $h\in\R$ $$|f(x, u)-f(x,0)|\le
\overline C |a(x)||u|+ \overline C|u|^{\overline \rho+1},
\leqno \dff{160405-1328}$$ $$|f(x, u+h)-f(x,u)|\le
\overline C |a(x)||h|+ \overline C\max(1,2^{\overline
\rho-1})(|u|^{\overline \rho}+ |h|^{\overline
\rho})|h|,\leqno\dff{290505-1204}$$
 $$| F(x,u)|\le
\overline C(|a(x)||u|^2/2+ |u|^{\overline
\rho+2}/(\overline \rho+2))+|u||
f(x,0)|,\leqno\dff{071105-1052}$$
$$\aligned
 &|F(x,u+h)-F(x,u)|\le\\& (|f(x,0)|+\overline C|a(x)|(|u|+|h|)+
 \overline C\max(1,2^{\overline \rho})(|u|^{\overline\rho+1}+|h|^{\overline \rho+1}))|h|,\endaligned\leqno\dff{060306-1329}$$and $$\aligned&| F(x,u+h)-
F(x,u)- f(x,u)h|\\&\le\bigl( \overline C|a(x)|+ \overline
C\max(1,2^{\overline \rho-1})(|u|^{\overline \rho}+
|h|^{\overline \rho})\bigr) |h|^2.
\endaligned\leqno\dff{081105-1936}$$Moreover, for every
measurable function $v\co \Omega\to\R$ both $\hat f(v)$ and
$\hat F(v)$ are measurable and for all measurable functions
$u$, $h\co\Omega\to\R$ $$|\hat f(u)|_{L^2}\le |\hat
f(0)|_{L^2}+\overline C(|a u|_{L^2}+ |u|^{\overline \rho
+1}_{L^{2(\overline \rho+1)}}), \leqno\dff{100405-1159}$$
$$\aligned|\hat f(u+h)&-\hat f(u)|_{L^2}\\&\le \overline
C|a h|_{L^2}+ \overline C\max(1,2^{\overline
\rho-1})(|u|^{\overline \rho}_
{L^{2(\overline\rho+1)}}+|h|^{\overline
\rho}_{L^{2(\overline\rho+1)}})
|h|_{L^{2(\overline\rho+1)}},
\endaligned\leqno\dff{100405-1200}
$$
 $$|\hat F(u)|_{L^1}\le
\overline C(\bigl|a |u|^2\bigr|_{L^1}/2+ |u|^{\overline
\rho+2}_{L^{\overline \rho +2}}/(\overline
\rho+2))+|u|_{L^2}|\hat
f(0)|_{L^2},\leqno\dff{081105-1947}$$
$$\aligned
 &|\hat F(u+h)-\hat F(u)|_{L^1}\le\\& (|\hat f(0)|_{L^2}+\overline C(|au|_{L^2}+|ah|_{L^2})+
 \overline C\max(1,2^{\overline \rho})(|u|_{L^{2(\overline\rho+1)}}^{\overline\rho+1}+ |h|_{L^{2(\overline\rho+1)}}^{\overline \rho+1}))|h|_{L^2},\endaligned\leqno\dff{060306-1352}$$
and $$\aligned&|\hat
F(u+h)-\hat F(u)-\hat f(u)h|_{L^1}\\&\le\bigl( \overline
C|a h|_{L^2}+ \overline C\max(1,2^{\overline
\rho-1})(|u|^{\overline \rho}_{L^{2(\overline\rho+1)}}+
|h|^{\overline
\rho}_{L^{2(\overline\rho+1)}})|h|_{L^{2(\overline\rho+1)}}\bigr)
|h|_{L^2}. \endaligned\leqno\dff{100405-1255}$$
Finally, if $\overline\rho$ is critical, then for every $r\in\ro N,\infty..$ there is a constant $C(r)\in\ro0,\infty..$
such that whenever $a=a_1+a_2$ with $a_1\in L^r(\Omega)$ and $a_2\in L^\infty(\Omega)$, then for all $u$, $h\in H^1_0(\Omega)$
$$\aligned|\hat f(u+h)-\hat f(u)|_{H^{-1}}&\le
C(r)(|a_1|_{L^r}+|a_2|_{L^\infty}) |h|_{L^{2}}\\& + C(r)(|u|^{\overline \rho}_
{L^{2^*}}+|h|^{\overline
\rho}_{L^{2^*}})
|h|_{L^{2}}.
\endaligned
\leqno\dff{270206-1552}$$
\endproclaim \proof For a.e.  $x\in \Om$ and all $u$, $h\in
\R$
  we have
$$f(x, u+h)-f(x,u)=\int_0^1\partial_u f(x,u+\theta h)h\,\D
\theta.$$ $$\aligned F(x,
u+h)-F(x,u)-f(x,u)h&=\int_0^1[f(x,u+\theta h)- f(x,u)]h
\,\D \theta\\& =\int_0^1[\int_0^1 \partial_uf(x,u+r\theta
h)\theta h\,\D r]h\,\D \theta \endaligned$$ These
equalities and the definition of $F$ imply
estimates~\rf{160405-1328}, ~\rf{290505-1204},
~\rf{071105-1052},~\rf{060306-1329} and~\rf{081105-1936}. Now well known arguments and H\"older inequality
yields the remaining assertions of  the proposition except~\rf{270206-1552}.
To prove~\rf{270206-1552}, let $r\in\ro N,\infty..$ and $u$, $h$ and $v\in H^1_0(\Omega)$ be arbitrary.
Then
$$
 \aligned
  &|\langle \hat f(u+h)-\hat f(u),v\rangle|\le \int_\Omega |f(x,(u+h)(x))-f(x,u(x))|\,|v(x)|\,\D x\\&+
  \overline C\int_\Omega |(ah)(x)|\,|v(x)|\,\D x+
  \overline C\max(1,2^{\overline\rho-1})\int_\Omega (|u(x)|^{\overline\rho}+|h(x)|^{\overline\rho})|h(x)||v(x)|\,\D x.
 \endaligned
$$
Now
$$
 \int_\Omega |(a_1h)(x)|\,|v(x)|\,\D x\le |a_1|_{L^r}|h|_{L^2}|v|_{L^{2r/(r-2)}}
$$
and
$$
 \int_\Omega |(a_2h)(x)|\,|v(x)|\,\D x\le |a_2|_{L^\infty}|h|_{L^2}|v|_{L^{2}}.
$$
Moreover, since $(1/2^*)+(1/2)+(1/N)=1$ and $N\overline \rho=2^*$ we also have
$$\int_\Omega (|u(x)|^{\overline\rho}+|h(x)|^{\overline\rho})|h(x)||v(x)|\,\D x\le(|u|_{L^{2^*}}^{\overline \rho}+|h|_{L^{2^*}}^{\overline \rho})|h|_{L^2}|v|_{L^{2^*}}
$$
Noting that $2r/(r-2)\le 2^*$ let $C$ be a common bound of the imbeddings $H^1_0(\Omega)\to L^s(\Omega)$ for $s\in\{2,2^*,2r/(r-2)\}$. Then we conclude
$$
 \aligned
  &|\langle \hat f(u+h)-\hat f(u),v\rangle|\\&\le
  \overline C C(|a_1|_{L^r}+|a_2|_{L^\infty})|h|_{L^2}|v|_{H^1}+
  \overline C\max(1,2^{\overline\rho-1})C(|u|_{L^{2^*}}^{\overline \rho}+|h|_{L^{2^*}}^{\overline \rho})|h|_{L^2}|v|_{H^1}.
 \endaligned
$$
Since
$$|\hat f(u+h)+\hat f(u)|_{H^{-1}}=\sup_{v\in H^1_0(\Omega)}
|\langle \hat f(u+h)-\hat f(u),v\rangle|$$
estimate~\rf{270206-1552} follows.\eproof

\proclaim{Standing Assumption}  For the rest of this paper,
we assume
 Hypothesis~\rf{221105-1756} and fix an $\eps\in\oi0,\infty..$. Let
$Z=H^1_0(\Omega)\times L^2(\Omega)$ and
$B=B_{\alpha,\beta,\eps}$ be defined as in
 Proposition~\rf{070505-1208}. Moreover, let $Y=L^2(\Omega)\times H^{-1}(\Omega)$.\endproclaim
\proclaim{Proposition~\dft{280505-0731}}
 Let $\overline C$, $\overline\rho\in \ro0,\infty..$
 and $a\co \Omega\to \R$ be a measurable function such
 that the assignments $u\mapsto |a| u$ and
 $u\mapsto |a|^{1/2} u$ induce bounded linear
 operators from $H^1_0(\Omega)$ to $L^2(\Omega)$.
 Suppose the function $f$ satisfies a
 $C^1$-Carath\'eodory condition and $\partial_u f$
 satisfies a $(\overline C,\overline \rho,a)$-growth
 condition. Moreover, suppose $f(\cdot,0)\in
 L^2(\Om)$. If $N\ge 3$, then assume also that
 $\overline\rho\le (2^*/2)-1$. Under these
 hypotheses, $ f$ induces a map
$\hat f\co H^1_0(\Om)\to L^2(\Om)$ which is  Lipschitzian
on bounded subsets of $H^1_0(\Om)$. The canonical primitive
$F$ of $f$ induces a map $\hat F\co H^1_0(\Omega)\to
L^1(\Omega)$. This map is Fr\'echet differentiable and
$D\hat F(u)[h]=\hat f(u)\cdot h$ for $u$ and $h\in
H^1_0(\Omega)$. The map $\Phi_f\co Z\to Z$,
$$\Phi_f(z)=(0,(1/\eps)\hat f(z_1)),\quad z=(z_1,z_2)\in
Z,\leqno\dff{230505-1458}$$ is bounded and Lipschitzian on
bounded subsets of $H^1_0(\Om)$. By $\pi_f$ we denote the
local semiflow $\Pi_{B,\Phi}$ on $Z$, where $\Phi=\Phi_f$.
This local semiflow does not explode in bounded subsets of
$H^1_0(\Om)$.\endproclaim \proof This follows from
Proposition~\rf{100405-0808}, the Sobolew imbedding
theorem, Proposition~\rf{070505-1208} and
Proposition~\rf{230505-1255}.\eproof
\remark{Remark~\dft{161105-1206}} By Lemma~\rf{161105-1152}
the hypotheses on the function $a$ imposed in
Proposition~\rf{280505-0731} are satisfied e.g. if $\tilde
a\in L^p_{\roman u}(\R^N)$ with $p\ge N$. \endremark
\remark{Remark~\dft{171105-0935}} The local semiflow
$\pi_{f}$ defined in Proposition~\rf{280505-0731} is, by
definition, the {\it local semiflow generated by solutions
of the damped
wave equation~\rf{261105-0951}\/}.%
\endremark

\head\dfs{100405-1955} Tail estimates and the existence of
attractors \endhead

 \proclaim{Proposition~\dft{l:diff}}
 Let $\overline\gamma\co\R^N\to\ci0,1..$ be a $C^1$-function such that $\sup_{x\in\R^N}
 (|\overline\gamma(x)|^2+|\nabla \overline\gamma(x)|^2)<\infty$. Set $\gamma=\overline\gamma^2$.
 Assume the hypotheses and notations of Proposition~\rf{280505-0731}.
 Fix $\delta\in \oi0,\infty..$,  and define the functions
$V=V_\gamma\co Z\to \R$ and $ V^*= V^*_\gamma\co Z\to \R$
by $$V(z)=(1/2)\int_{\Om}\gamma(x)\Psi_z(x)\,\D x $$ and $$
V^*(z)= \int_{\Om}\gamma(x)F(x,z_1(x))\,\D x$$ for
$z=(z_1,z_2)\in Z$. Here, for $z\in Z$ and $x\in\Omega$,
$$\Psi_z(x)=\eps|\delta z_1(x)+z_2(x)|^2+ (A\nabla
z_1)(x)\cdot \nabla z_1(x)+(\beta(x)-\delta
\alpha(x)+\delta^2\eps)|z_1(x)|^2.$$Let $\tau_0\in\oi
0,\infty..$, $I=\ci0,\tau_0..$ and $z\co I\to Z$ be a
solution of  $\pi_f$.
 Then the functions
$V\circ z$ and $ V^*\circ z$ are  differentiable and, for
$t\in I$, $$\aligned&(V\circ
z)'(t)=\\&\int_{\Om}\gamma(x)\bigl(\eps(\delta z_1+z_2)
(\delta z_2+(-(1/\eps)\alpha(x) z_2+(1/\eps)f(x,z_1(t)(x)
))\\ &+(-\delta \alpha(x)+\delta^2\eps)
z_1z_2-\delta\beta(x)z_1z_1\bigr)\,\D
x-\delta\int_{\Om}\gamma(x)(A\nabla( z_1))\cdot\nabla
 z_1\,\D x
 \\&-\int_{\Om}(\delta z_1+z_2)(A\nabla \gamma)\cdot\nabla
 z_1\,\D x\endaligned$$
$$( V^*\circ
z)'(t)=\int_{\Om}\gamma(x)f(x,z_1(t)(x))z_2(t)(x)\,\D x.$$
$$\aligned(V\circ z)'(t)+2\delta (V\circ
z)(t)&=\int_{\Om}\gamma(x)(2\delta\eps-\alpha(x))(\delta
z_1+z_2)^2\,\D x\\&+ \int_{\Om}\gamma(x)(\delta
z_1+z_2)f(x,z_1(t)(x))\,\D x \\&-\int_{\Om}(\delta
z_1+z_2)(A\nabla\gamma)\cdot \nabla z_1\,\D
x\endaligned\leqno\dff{280305-1000}$$\endproclaim \proof By
Proposition~\rf{100405-0808} we have that $V$ and $ V^*$
are defined and Fr\'echet differentiable on $Z$ and for all
$z=(z_1,z_2)$ and $\xi=(\xi_1,\xi_2)$ in $Z$ $$\multline
DV(z)[\xi]=\int_{\Om}\gamma(x)\bigl(\eps(\delta
z_1(x)+z_2(x)) (\delta \xi_1(x)+\xi_2(x))\\+ (A(x)\nabla
z_1(x))\cdot \nabla \xi_1(x) + (\beta(x)-\delta
\alpha(x)+\delta^2\eps)z_1(x)\xi_1(x)\bigr)\,\D x
\endmultline $$ and $$D
V^*(z)[\xi]=\int_{\Om}\gamma(x)f(x,z_1(x))\xi_1(x)\,\D x.$$
In particular, for $z=(z_1,z_2)\in D(B)$ and
$w=(w_1,w_2)\in Z$ we obtain, omitting the argument
$x\in\Om$ in some of the expressions below, $$\multline
DV(z)[B z+w]\\=\int_{\Om}\gamma(x)\bigl(\eps(\delta
z_1+z_2) (\delta (z_2+w_1)+(-(1/\eps)\alpha(x)
z_2-(1/\eps)\beta(x) z_1+(1/\eps)\Li z_1+w_2 ))\\+ (A\nabla
z_1)\cdot \nabla (z_2+w_1) +(\beta(x)-\delta
\alpha(x)+\delta^2\eps) z_1(z_2+w_1)\bigr)\,\D
x\endmultline $$ and $$D V^*(z)[B
z+w]=\int_{\Om}\gamma(x)f(x,z_1(x))(z_2+w_1)\,\D x.$$
Evaluating further we see that $$\multline DV(z)[B
z+w]\\=\int_{\Om}\gamma(x)\bigl(\eps(\delta z_1+z_2)
(\delta (z_2+w_1)+(-(1/\eps)\alpha(x) z_2-(1/\eps)\beta(x)
z_1+w_2 ))\\+ (A\nabla z_1)\cdot \nabla w_1
+(\beta(x)-\delta \alpha(x)+\delta^2\eps)
z_1(z_2+w_1)\bigr)\,\D x\\+
\int_{\Om}\gamma(x)\bigl((\delta z_1+z_2)\Li
 z_1+(A\nabla z_1)\cdot\nabla z_2\bigr)\,\D x.\endmultline $$
By Green's formula
$$\multline\int_{\Om}\gamma(x)\bigl((\delta z_1+z_2)\Li
 z_1+(A\nabla z_1)\cdot\nabla z_2\bigr)\,\D
 x=-\int_{\Om}\gamma(x)(A\nabla(\delta z_1+z_2))\cdot\nabla
 z_1\,\D x
 \\-\int_{\Om}(\delta z_1+z_2)(A\nabla \gamma)\cdot\nabla
 z_1\,\D x +\int_{\Om}\gamma(x)(A\nabla z_1)\cdot\nabla z_2\,\D
 x\\=-\int_{\Om}\gamma(x)(A\nabla(\delta z_1))\cdot\nabla
 z_1\,\D x
 -\int_{\Om}(\delta z_1+z_2)(A\nabla \gamma)\cdot\nabla
 z_1\,\D x\endmultline$$
 so we  obtain
 $$\multline DV(z)[B z+w]\\=\int_{\Om}\gamma(x)\bigl(\eps(\delta z_1+z_2)
(\delta (z_2+w_1)+(-(1/\eps)\alpha(x) z_2+w_2 ))\\+
(A\nabla z_1)\cdot \nabla w_1 +(-\delta
\alpha(x)+\delta^2\eps) z_1(z_2+w_1)+\beta(x)(z_1w_1-\delta
z_1z_1 )\bigr)\,\D x\\-\int_{\Om}\gamma(x)(A\nabla(\delta
z_1))\cdot\nabla
 z_1\,\D x
 -\int_{\Om}(\delta z_1+z_2)(A\nabla \gamma)\cdot\nabla
 z_1\,\D x.\endmultline $$
 Define the maps $W\co Z\times Z\to \R$ and $ W^*\co Z\times Z\to \R$  by
$$\multline W(z,w)=\int_{\Om}\gamma(x)\bigl(\eps(\delta
z_1+z_2) (\delta (z_2+w_1)+(-(1/\eps)\alpha(x) z_2+w_2
))\\+ (A\nabla z_1)\cdot \nabla w_1 +(-\delta
\alpha(x)+\delta^2\eps) z_1(z_2+w_1)+\beta(x)(z_1w_1-\delta
z_1z_1 )\bigr)\,\D x\\-\int_{\Om}\gamma(x)(A\nabla(\delta
z_1))\cdot\nabla
 z_1\,\D x
 -\int_{\Om}(\delta z_1+z_2)(A\nabla \gamma)\cdot\nabla
 z_1\,\D x\endmultline $$
 and
 $$ W^*(z,w)= \int_{\Om}\gamma(x)f(x,z_1(x))(z_2(x)+w_1(x))\,\D x$$
 for $(z,w)\in Z\times Z$.
In the particular case where $w_1=0$ we thus obtain
  $$\aligned W(z,w)&=\int_{\Om}\gamma(x)\bigl(\eps(\delta z_1+z_2)
(\delta z_2+(-(1/\eps)\alpha(x) z_2+w_2 ))\\
&+(-\delta \alpha(x)+\delta^2\eps)
z_1z_2-\delta\beta(x)z_1z_1\bigr)\,\D
x\\&-\int_{\Om}\gamma(x)(A\nabla(\delta z_1))\cdot\nabla
 z_1\,\D x
 -\int_{\Om}(\delta z_1+z_2)(A\nabla \gamma)\cdot\nabla
 z_1\,\D x\endaligned\leqno \dff{270305-1941}$$
and $$ W^*(z,w)=\int_{\Om}\gamma(x)f(x,z_1(x))z_2(x)\,\D
x.\leqno\dff{280305-0840}$$ Using
Hypothesis~\rf{221105-1756} and
Proposition~\rf{100405-0808} we see that $W$ and $ W^*$ are
continuous from $Z\times Z$ to $\R$ and so
Theorem~\rf{p:diff}, formulas~\rf{270305-1941}
and~\rf{280305-0840} and a straightforward computation
complete the proof.\eproof Consider the following
hypothesis. \proclaim{Hypothesis~\dft{231105-0903}}
\roster\item $\alpha_0>0$;\item
 $\overline C$, $\overline\rho$, $\overline
\tau\in\ro0,\infty..$ and $\overline\mu\in \oi0,\infty..$
are constants and $c\co\Om\to \ro0,\infty..$ is a function
with $c\in L^1(\Om)$. If $N\ge 3$, then  $\overline
\rho\le (2^*/2)-1$;\item
 $a\co \Omega\to \R$ is a measurable function
  such that the assignments $u\mapsto |a|
 u$ and $u\mapsto |a|^{1/2} u$ induce
 bounded linear operators from $H^1_0(\Omega)$
 to $L^2(\Omega)$; \item $f\co \Omega\times
 \R\to \R$  satisfies a
$C^1$-Carath\'eodory condition;\item
$F$ is the canonical primitive of $f$;
\item $\partial_u f$ satisfies a $(\overline C,\overline
\rho,a)$-growth condition;\item $| f(\cdot,0)|_{L^2}\le
\overline \tau$; \item $f(x,u)u-\overline\mu F(x,u)\le
c(x)$ and $F(x,u)\le c(x)$ for a.e.  $x\in \Om$ and every
$u\in\R$. \endroster  \endproclaim

A sufficient condition for the dissipativity
assumption~$(8)$ to hold is contained in the following
lemma: \proclaim{Lemma~\dft{121205-1154}} Let
$f\co\Omega\times\R\to\R$ satisfy a $C^0$-Carath\'eodory
condition and $F$ be the canonical primitive of $f$. Let
$\nu$, $\gamma\in\oi1,\infty..$ be constants and  $D\in
L^1(\Omega)$ be a function with $D(x)>0$ for all
$x\in\Omega$ and such that $F(x,u)\le D(x)$ for all
$x\in\Omega$ and all $u\in\R$. Assume also that the
function $u\mapsto (\gamma D(x)-F(x,u))^\nu$ is convex for
a.e. $x\in\Omega$.

Then  $f(x,u)u-\overline\mu F(x,u)\le c(x)$ and $F(x,u)\le
c(x)$ for a.e.  $x\in \Om$ and every $u\in\R$. Here,
$\overline \mu:=(1/\nu)$ and
$c(x):=\max(1,\gamma^\nu(\gamma-1)^{1-\nu}\nu^{-1}) D(x)$,
$x\in\Omega$. \endproclaim%

\proof Define $G(x,u)=-(\gamma D(x)-F(x,u))^\nu$ for
$x\in\Omega$ and $u\in\R$. Our convexity assumption implies
that the function $u\mapsto \partial_u G(x,u)$ is
nonincreasing and continuous for a.e. $x\in\Omega$. Notice
that whenever $h\co\R\to\R$ is continuous and nonincreasing
then $h(u)u\le\int_0^u h(s)\,\D s$ for all $u\in\R$. It
follows that, for a.e. $x\in\Omega$ and every $u\in\Omega$,
$$\aligned\nu f(x,u)(\gamma D(x)-F(x,u))^{\nu-1}u&\le
G(x,u)-G(x,0)\\&=-(\gamma D(x)-F(x,u))^\nu+(\gamma
D(x))^\nu.\endaligned\leqno\dff{191205-1052}$$ Since, for
a.e. $x\in\Omega$ and every $u\in\Omega$, $$\gamma
D(x)-F(x,u)\ge (\gamma-1)D(x)>0$$ we obtain
from~\rf{191205-1052} that%
$$\aligned\nu f(x,u)u&\le -(\gamma D(x)-F(x,u))+ (\gamma
D(x))^\nu((\gamma-1)D(x))^{1-\nu}\\&\le
F(x,u)+\gamma^\nu(\gamma-1)^{1-\nu}D(x).\endaligned$$
The lemma is proved.\eproof%
  Fix a $C^\infty$-function $\overline\vartheta\co
\R\to \ci 0,1..$ with $\overline\vartheta (s)=0$ for $s\in
\lo-\infty,1..$ and $\overline\vartheta (s)=1$ for $s\in
\ro2,\infty..$. Let $$\vartheta:=\overline\vartheta^2.$$
For $k\in\N$ let the functions
$\overline\vartheta_k\co\R^N\to\R$ and $\vartheta_k\co
\R^N\to\R$ be defined by
$$\overline\vartheta_k(x)=\overline\vartheta(|x|^2/k^2)\text{
and }\vartheta_k(x)=\vartheta(|x|^2/k^2),\quad x\in \R^N.$$
\proclaim{Theorem~\dft{290505-2035}} Assume
Hypothesis~\rf{231105-0903}. Choose $\delta$ and
$\nu\in\oi0,\infty..$ with $$\text{
$\nu\le\min(1,\overline\mu/2)$,
$\lambda_1-\delta\alpha_1>0$ and $\alpha_0-2\delta\eps\ge
0$.}\leqno\dff{280305-1006}$$
  Under these hypotheses, there is a constant $ c' \in\ro0,\infty..$ and for
every $R\in \ro0,\infty..$  there are constants $ M'=
M'(R)$,  $c_k=c_k(R)\in
\ro0,\infty..$, $k\in\N$ with $c_k\to 0$ for $k\to \infty$
and such that for every
  $\tau_0\in
\ro0,\infty..$  and every  solution $z(\cdot)$  of $\pi_f$
on $I=\ci0,\tau_0..$ with $|z(0)|_Z\le R$
$$\aligned&\int_{\Om}\bigl((\eps/2)|z_2(t)(x)|^2+
(A(x)\nabla z_1(t)(x))\cdot\nabla
z_1(t)(x)+(\beta(x)-\delta \alpha(x)
)|z_1(t)(x)|^2\bigr)\,\D x\\&\le  c' + M'e^{-2\delta \nu
t},\quad t\in I.\endaligned\leqno\dff{211105-1931}$$ If
$|z(t)|_Z\le R$ for $t\in I$, then
$$\aligned&\int_{\Om}\vartheta_k(x)\bigl((\eps/2)|z_2(t)(x)|^2+
(A(x)\nabla z_1(t)(x))\cdot\nabla
z_1(t)(x)\\&+(\beta(x)-\delta \alpha(x)
)|z_1(t)(x)|^2\bigr)\,\D x\\&\le c_k + M'e^{-2\delta \nu
t},\quad k\in \N,\, t\in
I.\endaligned\leqno\dff{010605-0834}$$ \endproclaim
\proclaim{Lemma~\dft{221105-1015}} Assume the hypotheses of
Theorem~\rf{290505-2035}. Let $\overline \gamma$, $\gamma$,
$V=V_\gamma$ and $ V^*= V^*_\gamma$ be as in
Proposition~\rf{l:diff}. For all $z\in Z$ and $x\in \Omega$
define $s(z)(x)=s_{\overline\gamma}(z)(x)$ by $$s(z)(x)=-2
\overline\gamma(x)z_1(x)(A(x)\nabla
\overline\gamma(x))\cdot\nabla z_1(x)-|z_1(x)|^2(A(x)\nabla
\overline\gamma(x))\cdot \overline\gamma(x).$$ Given
$\tau_0\in\ro0,\infty..$ and a solution $z(\cdot)$ of
$\pi_f$ on $I=\ci0,\tau_0..$, define
$$\eta(t)=\eta_\gamma(t)=V_\gamma(z(t))-
V^*_\gamma(z(t)),\, t\in I.$$ Then
$$\aligned&\eta'(t)+2\delta\nu\eta(t)\le2
\delta(\overline\mu-\nu)\int_{\Om}\gamma(x)c(x)\,\D x
\\&-\int_{\Om}(\delta z_1+z_2)(A\nabla\gamma)\cdot \nabla
z_1\,\D x-\delta(1-\nu)\int_\Omega s_{\overline
\gamma}(z(t))(x)\,\D x,\,t\in
I.\endaligned\leqno\dff{280305-1352}$$\endproclaim \proof
It is clear that, for all $z\in Z$ and $x\in\Omega$
$$|\overline\gamma(x)|^2(A(x)\nabla z_1(x))\cdot \nabla
z_1(x) =(A(x)\nabla(\overline \gamma z_1)(x))\cdot\nabla
(\overline \gamma z_1)(x) +s(z)(x).$$ Thus, by the
definition of $V$, $$\aligned
2&V(z)\ge\int_{\Om}\gamma(x)\bigl( (A(x)\nabla z_1(x))\cdot
\nabla z_1(x)+(\beta(x)-\delta
\alpha(x))|z_1(x)|^2\bigr)\,\D x\\&=\int_{\Om}\bigl(
(A(x)\nabla(\overline \gamma z_1)(x))\cdot\nabla (\overline
\gamma z_1)(x)+(\beta(x)-\delta
\alpha(x))|\overline\gamma(x)z_1(x)|^2\bigr)\,\D
x\\&+\int_\Omega s_{\overline\gamma}(z)(x)\,\D x\ge
(\lambda_1-\delta\alpha_1)|\overline\gamma z_1|^2_{L^2}+
\int_\Omega  s_{\overline\gamma}(z)(x)\,\D x\ge\int_\Omega
s_{\overline\gamma}(z)(x)\,\D x.\endaligned $$ Hence
$$\aligned(V\circ z)'(t)&+2\delta\nu (V\circ z)(t)=(V\circ
z)'(t)+2\delta (V\circ z)(t)-\delta(1-\nu)2V(z(t))\\&\le
(V\circ z)'(t)+2\delta (V\circ
z)(t)-\delta(1-\nu)\int_\Omega
s_{\overline\gamma}(z(t))(x)\,\D x.\endaligned$$ It follows
that $$\aligned &(V\circ z)'(t)+2\delta\nu (V\circ z)(t)+
\delta(1-\nu)\int_\Omega s_{\overline\gamma}(z(t))(x)\,\D
x\\&\le (V\circ z)'(t)+2\delta (V\circ
z)(t)\\&\le(2\delta\eps-\alpha_0)\int_{\Om}\gamma(x)(\delta
z_1+z_2)^2\,\D x+ \int_{\Om}\gamma(x)(\delta
z_1+z_2)f(x,z_1(t)(x))\,\D x \\&-\int_{\Om}(\delta
z_1+z_2)(A\nabla\gamma)\cdot \nabla z_1\,\D
x\le\delta\int_{\Om}\gamma(x) z_1f(x,z_1(t)(x))\,\D x\\&+
\int_{\Om}\gamma(x)z_2f(x,z_1(t)(x))\,\D x
-\int_{\Om}(\delta z_1+z_2)(A\nabla\gamma)\cdot \nabla
z_1\,\D x\\&\le
\delta\overline\mu\int_{\Om}\gamma(x)(F(x,z_1(t)(x))
+c(x))\,\D x -2\delta\nu( V^*\circ z)(t)\\&+2\delta\nu(
V^*\circ z)(t) +( V^*\circ z)'(t) -\int_{\Om}(\delta
z_1+z_2)(A\nabla\gamma)\cdot \nabla z_1\,\D
x=:S^*\endaligned$$
Now $$\aligned S^*&=\delta(\overline\mu-2\nu)\int_{\Om}\gamma(x)F(x,z_1(t)(x))\,\D
x+\delta\overline\mu\int_{\Om}\gamma(x)c(x)\,\D x
\\&+2\delta\nu( V^*\circ z)(t) +( V^*\circ z)'(t)
-\int_{\Om}(\delta z_1+z_2)(A\nabla\gamma)\cdot \nabla
z_1\,\D
x\\&\le\delta(2\overline\mu-2\nu)\int_{\Om}\gamma(x)c(x)\,\D
x\\&+2\delta\nu( V^*\circ z)(t) +( V^*\circ z)'(t)
-\int_{\Om}(\delta z_1+z_2)(A\nabla\gamma)\cdot \nabla
z_1\,\D x.\endaligned$$ This immediately
implies~\rf{280305-1352} and proves the lemma.\eproof
\demo{Proof of Theorem~\rf{290505-2035}} Let $\tau_0\in
\ro0,\infty..$ be arbitrary and  $z(\cdot)$ of $\pi_f$ be
an arbitrary solution of $\pi_f$ on $I=\ci0,\tau_0..$ with
$|z(0)|_Z\le R$. Set $\overline \gamma=\gamma\equiv1$. Then
$s_{\overline\gamma}(z(t))\equiv0$. Thus
Lemma~\rf{221105-1015} implies that $$\eta_\gamma'+2\delta
\nu \eta_\gamma\le \overline c.\leqno\dff{221105-1358}$$
where $\overline
c=2\delta(\overline\mu-\nu)\int_{\Om}\gamma(x)c(x)\,\D x$.
Differentiating the function
 $t\mapsto\eta_\gamma(t)e^{2\delta\nu t}$ and using~\rf{221105-1358} we obtain
 $$\eta_\gamma(t)\le (1/(2\delta\nu))\overline c[1-e^{-2\delta\nu
 t}]+\eta_\gamma(0)e^{-2\delta\nu t},\quad t\in
 I.\leqno\dff{221105-1401}$$
 Our assumptions imply that there is a continuous imbedding
$H^1_0(\Om)\to L^{\overline \rho+2}(\Om)$ with an imbedding
constant $C_2$. Let $L_\beta$, resp. $L_a$ be bounds on the
operators from $H^1_0(\Om)$ to $L^{2}(\Om)$ given by the
assignments $u\mapsto |\beta|^{1/2} u$, resp. $u\mapsto
|a|^{1/2} u$.
 Now a simple calculation using Proposition~\rf{100405-0808} shows that
 $$\aligned|\eta_\gamma(0)|&\le (1/2)(2\delta^2\eps R^2+2\eps R^2+ a_1 R^2+
 (L_\beta^2+\delta^2\eps)R^2)\\&+\overline C(L_a^2R^2/2+(C_2)^
 {\overline \rho+2}R^{\overline \rho+2}/(\overline \rho+2))+R\overline\tau
 =:\overline M.\endaligned\leqno\dff{160405-1348}$$
 The definitions of $V_\gamma$ and $ V^*_\gamma$ and our assumption on $F$ now
  imply that, for $t\in I$,
$$\aligned (1/2)&\int_{\Om}\bigl(\eps|\delta
z_1(t)(x)+z_2(t)(x)|^2+ (A(x)\nabla z_1(t)(x))\cdot\nabla
z_1(t)(x) \\&+(\beta(x)-\delta
\alpha(x)+\delta^2\eps)|z_1(t)(x)|^2\bigr)\,\D x\\&\le
(1/(2\delta\nu))\overline c[1-e^{-2\delta\nu
 t}]+\overline M e^{-2\delta\nu t}+ V^*_\gamma(z(t))\\&\le (1/(2\delta\nu))
 \overline c[1-e^{-2\delta\nu
 t}]+\overline M e^{-2\delta\nu t}+\int_\Omega c(x)\,\D x.\endaligned\leqno
 \dff{221105-1419}$$
Now, for $a_1$, $a_2\in \R$ we have
$$|a_1|^2=|(a_1+a_2)+(-a_2)|^2\le 2(|a_1+a_2|^2+|a_2|^2)$$
so $$ |a_1+a_2|^2\ge (1/2)|a_1|^2-|a_2|^2$$and thus setting
$a_1=z_2(t)(x)$ and $a_2=\delta z_1(t)(x)$
in~\rf{221105-1419} we obtain $$\aligned
(1/2)&\int_{\Om}\bigl((\eps/2)|z_2(t)(x)|^2+ (A(x)\nabla
z_1(t)(x))\cdot\nabla z_1(t)(x)\\&+(\beta(x)-\delta
\alpha(x))|z_1(t)(x)|^2\bigr)\,\D x\\&\le
(1/(2\delta\nu))\overline c[1-e^{-2\delta\nu
 t}]+\overline M e^{-2\delta\nu t}+\int_\Omega c(x)\,\D x.\endaligned
 \leqno\dff{221105-1420}$$
Setting $ c'=2((1/(2\delta\nu))\overline c+\int_\Omega
c(x)\,\D x)$ and $ M'=2\overline M$ we
obtain~\rf{211105-1931}.

Assume now that $|z(t)|_Z\le R$ for all $t\in I$. Let $k\in
\N$ be arbitrary and set $V_k=V_{\gamma_k}$, $ V^*_k=
V^*_{\gamma_k}$, $s_k(z)(x)=s_{\overline\gamma_k}(z)(x)$
and $\eta_k(t)=\eta_{\gamma_k}$, where
$\overline\gamma_k=\overline \vartheta_k$ and
$\gamma_k=\vartheta_k$. Since $\nabla
\vartheta_k(x)=(2/k^2)\vartheta'(|x|^2/k^2)x$ and $\nabla
\overline\vartheta_k(x)=(2/k^2)\overline\vartheta'(|x|^2/k^2)x$
we have $$\text{$\sup_{x\in\Omega}|\nabla \vartheta_k
(x)|\le C_\vartheta /k$ and $\sup_{x\in\Omega}|\nabla
\overline\vartheta_k (x)|\le C_{\overline\vartheta}
/k$}\leqno\dff{171105-0943}$$ where $C_\vartheta=2\sqrt
2\sup_{y\in\R}|\vartheta'(y)|$ and
$C_{\overline\vartheta}=2\sqrt
2\sup_{y\in\R}|\overline\vartheta'(y)|$.

  We thus obtain
 $$-\int_{\Om}(\delta z_1+z_2)(A\nabla\vartheta_k)\cdot \nabla z_1\,\D
 x\le
    a_1 (C_\vartheta/k)(\delta R+R)R
     \leqno\dff{160405-1236}$$
     and
$$-\delta(1-\nu)\int_\Omega s_k(z(t))(x)\,\D x\le a_1
\delta(1-\nu)
(2C_{\overline\vartheta}/k+C_{\overline\vartheta}^2/k^2)R^2.\leqno\dff{161105-1831}$$
 Set
 $$\aligned&\xi_k= 2\delta(\overline\mu-\nu)
\int_{\{x\in\Om\mid |x|\ge k\}}|c(x)|\,\D x\\&+
a_1(C_\vartheta/k)(\delta R+R)R+a_1\delta(1-\nu)
(2C_{\overline\vartheta}/k+C_{\overline\vartheta}^2/k^2)R^2.
\endaligned\leqno\dff{280305-1448}$$ Using
Lemma~\rf{221105-1015} we thus have that $$\eta_k'+2\delta
\nu \eta_k\le \xi_k,\quad k\in \N.\leqno\dff{100405-1905}$$
Differentiating the function
 $t\mapsto\eta_k(t)e^{2\delta\nu t}$ and using~\rf{100405-1905} we obtain
 $$\eta_k(t)\le (1/(2\delta\nu))\xi_k[1-e^{-2\delta\nu
 t}]+\eta_k(0)e^{-2\delta\nu t},\quad t\in
 I.\leqno\dff{100405-1911}$$
We have
 $$|\eta_k(0)|\le \overline M\leqno\dff{231105-1028}$$
 where $\overline M$ is as in~\rf{160405-1348}.
Using our assumptions on $\vartheta$ we obtain $$\aligned&
V^*(z(t))\le\int_{\Om}\vartheta_k(x)c(x)\,\D x \le
\int_{\{x\in\Om\mid |x|\ge k\}}c(x)\,\D x=:\zeta_k,\quad
t\in I. \endaligned\leqno\dff{280305-1710}$$ It follows
that, for $t\in I$,
$$\aligned&(1/2)\int_{\Om}\vartheta_k(x)\bigl(\eps|\delta
z_1(t)(x)+z_2(t)(x)|^2+ (A(x)\nabla z_1(t)(x))\cdot\nabla
z_1(t)(x)\\&+(\beta(x)-\delta \alpha(x)
+\delta^2\eps)|z_1(t)(x)|^2\bigr)\,\D
x\le(1/(2\delta\nu))\xi_k+ \overline Me^{-2\delta\nu
t}+\zeta_k.\endaligned\leqno\dff{160405-1915}  $$ As
before, this implies that
$$\aligned&(1/2)\int_{\Om}\vartheta_k(x)\bigl((\eps/2)|z_2(t)(x)|^2+
(A(x)\nabla z_1(t)(x))\cdot\nabla
z_1(t)(x)\\&+(\beta(x)-\delta \alpha(x)
)|z_1(t)(x)|^2\bigr)\,\D x\\&\le(1/(2\delta\nu))\xi_k+
\overline Me^{-2\delta\nu
t}+\zeta_k.\endaligned\leqno\dff{160405-1932}$$ Setting $
M'=2\overline M$ and $c_k=2((1/(2\delta\nu))\xi_k+
\zeta_k)$, $k\in \N$ we obtain~\rf{010605-0834}. The
theorem is proved. \eproof
\proclaim{Theorem~\dft{231105-1032}} Assume
Hypothesis~\rf{231105-0903}. Then $\pi_f$ is a global
semiflow. Moreover, there is a constant $C_{\pi_f}\in
\ro0,\infty..$ with the property that for every $z_0$ there
is a $t_{z_0}\in\ro0,\infty..$ such that $|z_0\pi_f t|_Z\le
C_{\pi_f} $ for all $t\in \ro t_{z_0},\infty..$.
Furthermore, every bounded subset of $Z$ is ultimately
bounded (rel. to $\pi_f$).

\endproclaim \proof Using the first part of
Theorem~\rf{290505-2035} together with
Lemma~\rf{111105-1838} (with $\kappa=\delta \alpha_1 $) we
conclude that for every $z_0\in Z$ there is a constant
$C_{z_0}\in\ro0,\infty..$ such that $|z_0\pi_f t|_Z\le
C_{z_0}$ for $t\in\ro0,\omega_{z_0}..$. Since $\pi_f$ does
not explode in bounded subsets of $Z$, this implies that
$\omega_{z_0}=\infty$, so $\pi_f$ is a global semiflow.
Similar arguments prove the other assertions of the theorem.
\eproof
Now consider the following alternative hypotheses:
\proclaim{Hypothesis~\dft{050306-2123}}
 $\overline\rho$ is subcritical and $\tilde a\in L^r_{\roman {loc}}(\R^N)$ for some $r\in \R$ with $r>\max(N,2)$.
\endproclaim
\proclaim{Hypothesis~\dft{050306-2125}}
 $\overline \rho$ is critical, $a\in L^r(\Omega)+L^\infty(\Omega)$ for some $r\in\ro N,\infty..$ and~\rf{050306-0725} is satisfied.
 \endproclaim
\proclaim{Lemma~\dft{050306-2033}}
 Let $\tilde N$
 be an arbitrary ultimately bounded set in $Z=H^1_0(\Omega)\times L^2(\Omega)$ (relative to $\pi_f$), $(z_n)_n$ be an
 arbitrary sequence in $\tilde N$  and $(t_n)_n$ be a
 sequence in $\ro0,\infty..$ with $t_n\to \infty$.
  \roster
  \item  if Hypothesis~\rf{050306-2123} holds, then the sequence $(z_n\pi_f t_n)_n$ has a subsequence which converges in $Z$.
  \item if Hypothesis~\rf{050306-2125} holds, then $(z_n\pi_f t_n)_n$ has a subsequence which converges in $Y=L^2(\Omega)\times H^{-1}(\Omega)$.
 \endroster
\endproclaim
\proof
There is
a $t_{\tilde N}$ and an $R\in \ro0,\infty..$ such that
$|z\pi_f t|_Z\le R$ for all $z\in \tilde N$ and all
$t\in\ro t_{\tilde N},\infty..$. We may assume that $t_n\ge
t_{\tilde N}$ and therefore, replacing $z_n$ by $z_n\pi_f
t_{\tilde N}$ and $t_n$ by $t_n-t_{\tilde N}$ we may assume
that
   $|z_n\pi_{f}t|_Z\le R$ for all $n\in \N$ and $t\in
\ci0,t_n..$. For $n\in\N$ and $t\in\ci0,t_n..$ let $u_n(t)$
be the first component of $z_n\pi_{f}t$. Let $\tau_0\in
\oi0,\infty..$ be arbitrary to be determined later.
 Then there an $n_0(\tau_0)\in \N$ such that
$t_n\ge 2\tau_0$ for all $n\in \N$ with $n\ge n_0(\tau_0)$.
For such $n$ we have $$\aligned &z_n\pi_{f}t_n=
T(\tau_0)z_n\pi_{f}(t_n-\tau_0)\\&+
\int_0^{\tau_0}T(\tau_0-s)(0,(1/\eps)(\hat
f(u_n(t_n-\tau_0+s))- \hat
f((1-\overline\vartheta_k)u_n(t_n-\tau_0+s))))\,\D s \\&+
\int_0^{\tau_0}T(\tau_0-s)(0,(1/\eps) \hat
f((1-\overline\vartheta_k) u_n(t_n-\tau_0+s))\,\D
s\endaligned$$ We have
$$|T(\tau_0)z_n\pi_{f}(t_n-\tau_0)|_Z\le Me^{-\mu \tau_0}R,
\quad n\ge n_0(\tau_0).\leqno\dff{070605-1856}$$ Since
$\sup_{k\in
\N}|\overline\vartheta_k|_{W^{1,\infty}(\R^N)}<\infty$ it
follows from Lemma~\rf{160505-0804} that
$$\sup_{k,n\in\N}\sup_{t\in\ci0,t_N..}(|u_n(t)|_{H^1_0}+ |(1-\overline\vartheta_k)u_n(t)
|_{H^1_0})<\infty.$$ It follows from our hypotheses and
from Proposition~\rf{280505-0731} that there is an
$L\in\oi0,\infty..$ such that for all $k\in\N$, $n\in\N$
and $t\in\ci0,t_n..$ $$|\hat f(u_n(t))- \hat
f((1-\overline\vartheta_k)u_n(t))|_{L^2}\le
L|\overline\vartheta_k u_n(t)|_{H^1_0}.$$ This implies that
$$\aligned|\int_0^{\tau_0}&T(\tau_0-s)(0,(1/\eps)(\hat
f(u_n(t_n-\tau_0+s))- \hat
f((1-\overline\vartheta_k)u_n(t_n-\tau_0+s))))\,\D
s|_Z\\&\le \sup_{s\in\ci0,\tau_0..}|\overline\vartheta_k
u_n(t_n-\tau_0+s)|_{H^1_0}(1/\eps)L M
\int_0^{\tau_0}e^{-\mu (\tau_0-s)}\,\D s\\& \le
(LM/(\mu\eps))\sup_{s\in\ci0,\tau_0..}|\overline\vartheta_k
u_n(t_n-\tau_0+s)|_{H^1_0},\quad n\ge
n_0(\tau_0).\endaligned\leqno\dff{070605-1857}$$ Now use
Lemma~\rf{111105-1838} with $\kappa=\delta\alpha_1$. Let
$c>0$ be as in that Lemma. It follows that, for $k$,
$n\in\N$ and $t\in \ci0,t_n..$ $$\aligned
c&|\overline\vartheta_k u_n(t)|_{H^1_0}^2\le \is
A\nabla(\overline \vartheta_k u_n(t)).. \nabla(\overline
\vartheta_k u_n(t))..+\langle \beta\overline \vartheta_k
u_n(t),\overline \vartheta_k
u_n(t)\rangle\\&-\delta\alpha_1\is \overline \vartheta_k
u_n(t)..\overline \vartheta_k u_n(t).. \\&\le\is
A\nabla(\overline \vartheta_k u_n(t)).. \nabla(\overline
\vartheta_k u_n(t))..+\langle \beta\overline \vartheta_k
u_n(t),\overline \vartheta_k u_n(t)\rangle\\&-\delta\is
\alpha\overline \vartheta_k u_n(t)..\overline \vartheta_k
u_n(t)..\\&= \int_{\Om}\vartheta_k(x)\bigl( \is A\nabla
u_n(t).. \nabla u_n(t)..+(\beta(x)-\delta \alpha(x)
)|u_n(t)(x)|^2\bigr)\,\D x\\&+2\is\overline \vartheta_k
A\nabla u_n(t).. u_n(t)\nabla\overline \vartheta_k ..+\is
u_n(t) A\nabla\overline \vartheta_k ..
u_n(t)\nabla\overline \vartheta_k..\\&\le c_k +
M'e^{-2\delta \nu
t}+a_1(2C_{\overline\vartheta}/k+C_{\overline\vartheta}^2/k^2)R^2
\endaligned\leqno\dff{080605-0848}$$  Now, if $n\ge
n_0(\tau_0)$ and $s\in \ci 0,\tau_0..$ then
$t=t_n-\tau_0+s\ge \tau_0$ so \rf{080605-0848} implies that
 $$\sup_{n\ge n_0(\tau_0)}\sup_{s\in\ci0,\tau_0..}|\overline\vartheta_k
  u_n(t_n-\tau_0+s)|_{H^1_0}\to 0$$for $k\to\infty$ and $\tau_0\to \infty$.
 It follows that the right hand sides of~\rf{070605-1856} and~\rf{070605-1857}
  can be made as small as we wish, by taking $k\in \N$ and $\tau_0>0$
   sufficiently large. Therefore, a standard argument using Kuratowski
    measure of noncompactness implies that the sequence $(z_n\pi_{f}t_n)_n$
     has a subsequence which converges in $Z$ (resp. in $Y$) provided we can prove that,
 for every $k\in \N$ and $\tau_0\in \oi0,\infty..$ the set
 $$K_0:=\{\,T(\tau_0-s)(0,(1/\eps)\hat f((1-\overline \vartheta_k)u_n(t_n-\tau_0+s))\mid n\ge
n_0(\tau_0),\, s\in\ci0,\tau_0..\,\}$$
 is relatively compact in $Z$ (resp. in $Y$).

  Let $(z_l)_l$ be a sequence in $K_0$. It follows that for every $l\in \N$ there are $n_l\in\N$ $s_l\in\ci0,\tau_0..$ with $z_l=T(\tau_0-s_l)(0,(1/\eps)\hat f(v_l))$
 where $v_l=(1-\overline \vartheta_k)u_{n_l}(t_{n_l}-\tau_0+s_l)$. By choosing subsequences if necessary we may assume that $s_l\to s_\infty$ for some $s_\infty\in \ci0,\tau_0..$. By Proposition~\rf{150505-1711}
 $(v_l)_l$  is compact in $L^s(\Omega)$ for each $s\in\ro2,\infty..$ such that $s\in\ro2,2^*..$ if $N\ge 3$.

 First suppose that Hypothesis~\rf{050306-2123} holds. Then $s\in\ro2,2^*..$ for $s\in \{2r/(r-2), 2(\overline\rho+1)\}$. Taking subsequences if necessary, we may thus assume that there is a $v\in H^1_0(\Omega)$  such that
 $(v_l)_l$ converges  to $v$ weakly in $H^1_0(\Omega)$ and strongly in $L^s(\Omega)$ for $s\in \{2r/(r-2), 2(\overline\rho+1)\}$.
 Moreover, whenever $x\in\Omega$ and $|x|\ge \sqrt 2 k$ then $v_l(x)=0$ for all $l\in\N$, and so we may assume that $v(x)=0$. Thus
 $$a(x)(v_l(x)-v(x))=a_1(x)(v_l(x)-v(x)), \quad l\in\N,\,x\in\Omega \leqno\dff{060306-1640}$$
 where $a_1\co\R^N\to \R$ is defined by $a_1(x)=\tilde a(x)$ if $x\in \overline\Omega$ and $|x|\le \sqrt 2 k$ and $a_1(x)=0$ otherwise. Note that $a_1\in L^r(\R^N)$ so the map $L^{2r/(r-2)}(\Omega)\to L^2(\Omega)$, $h\mapsto a_1 h$ is defined, linear and bounded.
 Now~\rf{100405-1200} and~\rf{060306-1640} imply that
 $|\hat f(v_l)-\hat f(v)|_{L^2}\to 0$ as $l\to \infty$.
 This clearly implies that $|z_l-T(\tau_0-s_\infty)(0,(1/\eps)\hat f(v))|_Z\to 0$ as $l\to \infty$.

 Now suppose Hypothesis~\rf{050306-2125}. Then $2\in\ro2,2^*..$ for $N\ge 3$.
 Taking subsequences if necessary, we may thus assume that there is a $v\in H^1_0(\Omega)$  such that
 $(v_l)_l$ converges  to $v$ weakly in $H^1_0(\Omega)$ and strongly in $L^2(\Omega)$. Using~\rf{270206-1552} we obtain that
 $|\hat f(v_l)-\hat f(v)|_{H^{-1}}\to 0$ as $l\to \infty$.
 Proposition~\rf{040306-2207} now implies that
 $|z_l-T(\tau_0-s_\infty)(0,(1/\eps)\hat f(v))|_Y\to 0$ as $l\to \infty$.
  The lemma
is proved.
\eproof
 We can now prove the first main result of this
paper. \proclaim{Theorem~\dft{231105-1037}} Assume
Hypotheses~\rf{231105-0903} and~\rf{050306-2123}. Then
$\pi_f$ is a global semiflow and it has a global attractor.
\endproclaim \proof This is an immediate consequence of
Theorem~\rf{231105-1032}, Lemma~\rf{050306-2033} and
Proposition~\rf{211105-1325}.\eproof
We will now treat the critical case.
\proclaim{Proposition~\dft{060306-1929}}
 Assume Hypotheses~\rf{231105-0903} and~\rf{050306-2125}.
 Let $C_{\dfc{070306-0708}}\in\ro0,\infty..$ be arbitrary.  Then there is a constant $C_{\dfc{070306-0709}}\in\ro0,\infty..$ such that whenever $t\in\ro0,\infty..$ and $z_1$ and $z_2\in Z$ are such that $|z_1|_{Z}\le C_{\rf{070306-0708}}$ and
 $|z_2|_{Z}\le C_{\rf{070306-0708}}$
 then
 $$|z_1\pi_f t-z_2\pi_f t|_{Y}\le C_{\rf{070306-0709}}e^{C_{\rf{070306-0709}}t}|z_1-z_2|_{Y}.$$

\endproclaim
\proof
 By Theorem~\rf{290505-2035} and Lemma~\rf{111105-1838} there is a constant $C_{\dfc{070306-0646}}\in\ro0,\infty..$ such that whenever $z\in Z$ and $|z|_Z\le C_{\rf{070306-0708}}$ then $|z\pi_f t|_{Z}\le C_{\rf{070306-0646}}$ for all $t\in\ro0,\infty..$.
 By~\rf{270206-1552} we now obtain a constant $C_{\dfc{070306-0647}}\in\ro0,\infty..$ such that $|\hat f(u_1)-\hat f(u_2)|_{H^{-1}}\le C_{\rf{070306-0647}}|u_1-u_2|_{L^2}$ for all $z_1=(u_1,v_1)$, $z_2=(u_2,v_2)\in Z$ with $|z_1|_{Z}\le C_{\rf{070306-0646}}$ and
 $|z_2|_{Z}\le C_{\rf{070306-0646}}$. Now Proposition~\rf{040306-2207}, the variation-of-constants formula and Gronwall's lemma complete the proof.
\eproof
\proclaim{Theorem~\dft{070306-0658}}
 Assume Hypotheses~\rf{231105-0903} and~\rf{050306-2125}. Then $\pi_f$ is asymptotically compact.
\endproclaim
\proof
 We use an ingenious method due to J. Ball, cf. \cite{\rfa{Ba}, \rfa{MRW}, \rfa{Ra}}.

 Let $\tilde N$ be a $\pi_f$-ultimately bounded subset of $Z$. Then there is a $t_{\tilde N}\in\ro0,\infty..$ and a $C_{\dfc{070306-0719}}\in\ro0,\infty..$  such that $|z\pi_f t|\le C_{\rf{070306-0719}}$ whenever $z\in \tilde N$ and $t\ge t_{\tilde N}$. Let $(z_n)_n$ be an arbitrary sequence in $\tilde N$ and $(t_n)_n$ be an arbitrary sequence in $\ro0,\infty..$ with $t_n\to \infty$ as $n\to \infty$. We must prove that a subsequence of $(z_n\pi_f t_n)_n$ converges strongly in $Z$. Now using Lemma~\rf{050306-2033} and Cantor's diagonal procedure we see that there is a strictly increasing sequence $(n_k)_k$ in $\N$ and for every $l\in\Z$ with $l\ge 0$ there  are a $k_0(l)\in\N$ and a $w_l\in Z$ with $|w_l|\le C_{\rf{070306-0719}}$  such that $t_{n_k}-l\ge t_{\tilde N}$ for $k\ge k_0(l)$ and the sequence
 $(z_{n_k}\pi_f (t_{n_k}-l))_{k\ge k_0(l)}$  converges to $w_l$ weakly in $Z$ and strongly in $Y$. By Proposition~\rf{060306-1929}, for every $l\in\N$ and $t\in\ro0,\infty..$,
 $$
  \text{$(z_{n_k}\pi_f (t_{n_k}-l))\pi_f t\to w_l\pi_f t$, as $k\to \infty$, strongly in $Y$.}
  \leqno\dff{070306-0905}
 $$ This shows that $w_l\pi_f l=w_0$
 for all $l\in \N$.
 Now define the function $\Cal F\co Z\to \R$ by
 $$
  \Cal F(z)=V(z)-V^*(z),\quad z\in Z
 $$
 where $V$ and $V^*$ are as in Proposition~\rf{l:diff} with $\gamma\equiv0$ and $\delta\in\oi0,\infty..$ such that $\lambda-\delta\alpha_1>0$ and $\alpha_0-2\delta\eps\ge0$.
 Using~\rf{081105-1947} we see that there is a constant  $C_{\dfc{070306-1111}}\in\ro0,\infty..$ such that
 $$
  \sup_{z\in Z, |z|_Z\le C_{\rf{070306-0719}}}|\Cal F (z)|\le C_{\rf{070306-1111}}.
 $$
 Note that $\Psi\co Z\to Z$, $(u,v)\mapsto (u,\delta u+v)$, is an isomorphism of normed spaces.
 Thus
 $$
  [(u_1,v_1),(u_2,v_2)]:= \eps\langle \delta u_1+v_1,\delta u_2+v_2\rangle+\langle A\nabla u_1,u_2\rangle+\langle(\beta-\delta\alpha+\delta^2\eps)u_1,u_2\rangle
 $$
 defines a scalar product on $Z$ whose norm $z\mapsto \|z\|:= \sqrt{[z,z]}$ is equivalent to the usual norm on $Z$.
 Note that $\Cal F(z)=\|z\|^2-V^*(z)$ for $z\in Z$.

 Let $\zeta=(\zeta_1,\zeta_2)\co\ro0,\infty..\to Z$  be an arbitrary solution of $\pi_f$.
 Proposition~\rf{l:diff} implies that  the function $\Cal F\circ \zeta$ is continuously differentiable and for every $t\in\ro0,\infty..$
 $$
  \aligned
  &(\Cal F\circ \zeta)'(t)+2\delta \Cal F(\zeta(t))=
  \int_{\Om}(2\delta\eps-\alpha(x))(\delta
  \zeta_1(t)(x)+\zeta_2(t)(x))^2\,\D x\\&+ \int_{\Om}\delta
  \zeta_1(t)(x)f(x,\zeta_1(t)(x))\,\D x -2\delta\int_{\Om}F(x,\zeta_1(t)(x))\,\D x.
  \endaligned
 $$
 It follows that for every $t\in \ro0,\infty..$
 $$
  \aligned
   &\Cal F(\zeta(t))=e^{-2\delta t}\Cal F(\zeta(0))\\&+\int_0^t e^{-2\delta (t-s)}\left(\int_{\Om}(2\delta\eps-\alpha(x))(\delta
  \zeta_1(s)(x)+\zeta_2(s)(x))^2\,\D x\right)\,\D s\\&+ \int_0^t e^{-2\delta (t-s)}\left(\int_{\Om}\delta
  \zeta_1(s)(x)f(x,\zeta_1(s)(x))\,\D x -2\delta\int_{\Om}F(x,\zeta_1(s)(x))\,\D x\right)\,\D s.
  \endaligned
  \leqno\dff{070306-0842}
 $$
 Fix $l\in\N$ and, for $k\ge k_0(l)$,
 let $\zeta_k(t)=(z_{n_k}\pi_f(t_{n_k}-l))\pi_f t$ and
  $\zeta(t)=w_l\pi_f t$ for $t\in\ro0,\infty..$. Then~\rf{070306-0842} with $t=l$ implies that
 $$
  \aligned
   &\|z_{n_k}\pi_f(t_{n_k})\|^2 -V^*(z_{n_k}\pi_f(t_{n_k}))=e^{-2\delta l}\Cal F(z_{n_k}\pi_f(t_{n_k}-l))\\&+\int_0^l e^{-2\delta (l-s)}\left(\int_{\Om}(2\delta\eps-\alpha(x))(\delta
  \zeta_{k,1}(s)(x)+\zeta_{k,2}(s)(x))^2\,\D x\right)\,\D s\\&+ \int_0^l e^{-2\delta (l-s)}\left(\int_{\Om}\delta
  \zeta_{k,1}(s)(x)f(x,\zeta_{k,1}(s)(x))\,\D x -2\delta\int_{\Om}F(x,\zeta_{k,1}(s)(x))\,\D x\right)\,\D s.
  \endaligned
  \leqno\dff{070306-0916}
 $$ and
 $$
  \aligned
   &\|w_0\|^2-V^*(w_0)=e^{-2\delta l}\Cal F(w_l)\\&+\int_0^l e^{-2\delta (l-s)}\left(\int_{\Om}(2\delta\eps-\alpha(x))(\delta
  \zeta_{1}(s)(x)+\zeta_{2}(s)(x))^2\,\D x\right)\,\D s\\&+ \int_0^l e^{-2\delta (l-s)}\left(\int_{\Om}\delta
  \zeta_{1}(s)(x)f(x,\zeta_{1}(s)(x))\,\D x -2\delta\int_{\Om}F(x,\zeta_{1}(s)(x))\,\D x\right)\,\D s.
  \endaligned
  \leqno\dff{070306-0926}
 $$
 Using~\rf{060306-1352} and~\rf{270206-1552} we see that
 $$V^*(z_{n_k}\pi_f(t_{n_k}))\to V^*(w_0)\leqno\dff{070306-1516}$$
 and
 $$
  \aligned
   &\int_0^l e^{-2\delta (l-s)}\left(\int_{\Om}\delta
  \zeta_{k,1}(s)(x)f(x,\zeta_{k,1}(s)(x))\,\D x -2\delta\int_{\Om}F(x,\zeta_{k,1}(s)(x))\,\D x\right)\,\D s\\&\to
  \int_0^l e^{-2\delta (l-s)}\left(\int_{\Om}\delta
  \zeta_{1}(s)(x)f(x,\zeta_{1}(s)(x))\,\D x -2\delta\int_{\Om}F(x,\zeta_{1}(s)(x))\,\D x\right)\,\D s
  \endaligned
  \leqno\dff{070306-1517}
 $$
 as $k\to\infty$.
 We claim that
 $$
  \aligned
   &\limsup_{k\to \infty}\int_0^l e^{-2\delta (l-s)}\left(\int_{\Om}(2\delta\eps-\alpha(x))(\delta
   \zeta_{k,1}(s)(x)+\zeta_{k,2}(s)(x))^2\,\D x\right)\,\D s\\&\le
   \int_0^l e^{-2\delta (l-s)}\left(\int_{\Om}(2\delta\eps-\alpha(x))(\delta
  \zeta_{1}(s)(x)+\zeta_{2}(s)(x))^2\,\D x\right)\,\D s.
  \endaligned
  \leqno\dff{070306-1126}
 $$
 In fact, since $\alpha(x)-2\delta \eps\ge0$ for all $x\in\Omega$
 we have by Fatou's lemma
 $$
  \aligned
   &\limsup_{k\to \infty}\int_0^l e^{-2\delta (l-s)}\left(\int_{\Om}(2\delta\eps-\alpha(x))(\delta
   \zeta_{k,1}(s)(x)+\zeta_{k,2}(s)(x))^2\,\D x\right)\,\D s\\&=
   -\liminf_{k\to \infty}\int_0^l e^{-2\delta (l-s)}\left(\int_{\Om}(\alpha(x)-2\delta\eps)(\delta
   \zeta_{k,1}(s)(x)+\zeta_{k,2}(s)(x))^2\,\D x\right)\,\D s\\&\le
   -\int_0^l e^{-2\delta (l-s)}\liminf_{k\to\infty}\left(\int_{\Om}(\alpha(x)-2\delta\eps)(\delta
   \zeta_{k,1}(s)(x)+\zeta_{k,2}(s)(x))^2\,\D x\right)\,\D s
  \endaligned
  \leqno\dff{070306-1524}
  $$
  Let $s\in\ci0,l..$ be arbitrary. Since $((\zeta_{k,1}(s), \zeta_{k,2}(s)))_k$ converges to $(\zeta_1(s),\zeta_2(s))$ weakly in $Z$ and $\Psi $ is continuous, linear, hence weakly continuous, it follows that $((\zeta_{k,1}(s), \delta \zeta_{k,1}(s)+\zeta_{k,2}(s)))_k$ converges to $(\zeta_1(s),\delta\zeta_1(s)+\zeta_2(s))$ weakly in $Z$. It follows that for every $v\in L^2(\Omega)$
  $$\langle v,\delta \zeta_{k,1}(s)+\zeta_{k,2}(s)\rangle\to
  \langle v,\delta \zeta_{1}(s)+\zeta_{2}(s)\rangle\text{ as $k\to\infty$.}$$
  Taking $v=(\alpha-2\delta\eps)(\delta\zeta_1(s)+\delta\zeta_2(s))$
  we thus obtain
  $$\aligned
  &|(\alpha-2\delta\eps)^{1/2}(\delta\zeta_1(s)+\delta\zeta_2(s))|_{L^2}^2\\&
  =\langle(\alpha-2\delta\eps)^{1/2}(\delta\zeta_1(s)+\delta\zeta_2(s)),
  (\alpha-2\delta\eps)^{1/2}(\delta\zeta_1(s)+\delta\zeta_2(s))\rangle\\&=
  \lim_{k\to\infty}\langle(\alpha-2\delta\eps)^{1/2} (\delta\zeta_1(s)+\delta\zeta_2(s)),
  (\alpha-2\delta\eps)^{1/2}(\delta\zeta_{k,1}(s)+\delta\zeta_{k,2}(s))\rangle\\&
  \le|(\alpha-2\delta\eps)^{1/2} (\delta\zeta_1(s)+\delta\zeta_2(s))|_{L^2}\liminf_{k\to\infty}| (\alpha-2\delta\eps)^{1/2}(\delta\zeta_{k,1}(s)+\delta\zeta_{k,2}(s))|_{L^2}
  \endaligned
   $$
  and so
  $$\aligned
   &\left(\int_{\Om}(\alpha(x)-2\delta\eps)(\delta
   \zeta_{1}(s)(x)+\zeta_{2}(s)(x))^2\,\D x\right)\\&\le
   \liminf_{k\to\infty}\left(\int_{\Om}(\alpha(x)-2\delta\eps)(\delta
   \zeta_{k,1}(s)(x)+\zeta_{k,2}(s)(x))^2\,\D x\right).
      \endaligned
      \leqno\dff{070306-1551}
  $$
  Inequalities~\rf{070306-1551} and~\rf{070306-1524} prove~\rf{070306-1126}.
 Using~\rf{070306-0916}, \rf{070306-0926}, \rf{070306-1516}, \rf{070306-1517} and~\rf{070306-1126} we obtain
 $$\aligned
  &\limsup_{k\to\infty}\|z_{n_k}\pi_f(t_{n_k})\|^2-V^*(w_l)\le
  e^{-2\delta l}C_{\rf{070306-1111}}\\&
  +\int_0^l e^{-2\delta (l-s)}\left(\int_{\Om}(2\delta\eps-\alpha(x))(\delta
  \zeta_{1}(s)(x)+\zeta_{2}(s)(x))^2\,\D x\right)\,\D s\\&+ \int_0^l e^{-2\delta (l-s)}\left(\int_{\Om}\delta
  \zeta_{1}(s)(x)f(x,\zeta_{1}(s)(x))\,\D x -2\delta\int_{\Om}F(x,\zeta_{1}(s)(x))\,\D x\right)\,\D s\\&=
  e^{-2\delta l}C_{\rf{070306-1111}} +\|w_0\|^2-V^*(w_0)-e^{-2\delta l}\Cal F(w_l)\le 2e^{-2\delta l}C_{\rf{070306-1111}}+\|w_0\|^2-V^*(w_0).
 \endaligned
 $$
 Thus for every $l\in\N$
 $$
 \limsup_{k\to\infty}\|z_{n_k}\pi_f(t_{n_k})\|^2\le 2e^{-2\delta l}C_{\rf{070306-1111}}+\|w_0\|^2
 $$
 so
 $$
  \limsup_{k\to\infty}\|z_{n_k}\pi_f(t_{n_k})\|\le\|w_0\|.
 $$
 Since $(z_{n_k}\pi_f(t_{n_k}))_k$ converges to $w_0$ weakly in $(Z,[\cdot,\cdot])$ we have
 $$
  \liminf_{k\to\infty}\|z_{n_k}\pi_f(t_{n_k})\|\ge\|w_0\|.
 $$
 Altogether we obtain
 $$\lim_{k\to\infty}\|z_{n_k}\pi_f(t_{n_k})\|=\|w_0\|.$$
 This implies that $(z_{n_k}\pi_f(t_{n_k}))_k$ converges to $w_0$ strongly in $Z$ and completes the proof.
\eproof
 We can now prove the second main result of this
paper. \proclaim{Theorem~\dft{070306-1631}} Assume
Hypotheses~\rf{231105-0903} and~\rf{050306-2125}. Then
$\pi_f$ is a global semiflow and it has a global attractor.
\endproclaim \proof This is an immediate consequence of
Theorem~\rf{231105-1032}, Theorem~\rf{070306-0658} and
Proposition~\rf{211105-1325}.\eproof

\Refs
\ref\no \dfa{AB1} \by W. Arendt and C. J. K. Batty
\paper Exponential stability of a diffusion equation with
absorption \jour Differential and Integral Equations \vol
6\yr 1993\pages 1009--1024
\endref
\ref\no \dfa{AB2} \bysame
\paper Absorption semigroups and Dirichlet boundary
conditions \jour Math. Ann. \vol 295\yr 1993\pages 427--448
\endref
\ref\no \dfa{ACDR} \by J. M. Arrieta, J. W. Cholewa,
T. D\l otko and A. Rodriguez-Bernal \paper Asymptotic
behavior and attractors for reaction diffusion equations in
unbounded domains \jour Nonlinear Analysis \vol 56\yr
2004\pages 515--554
\endref
\ref\no \dfa{BV} \by A. V.
Babin and M. I. Vishik \paper Regular attractors of
semigroups and evolution equations \jour J. Math. Pures
Appl. \vol 62\yr 1983\pages 441--491
\endref
\ref\no \dfa{Ba} \by J. M. Ball\paper Global attractors for damped semilinear wave equations. Partial differential equations and applications\jour  Discrete Contin. Dyn. Syst. \vol 10  \yr 2004  \pages 31--52
\endref
\ref\no \dfa{CH} \by T. Cazenave and A. Haraux \book
An Introduction to Semilinear Evolution Equations \publ
Clarendon Press \publaddr Oxford \yr 1998
\endref
\ref\no
\dfa{DCh} \by J. Cholewa and T. D\l otko  \book Global
Attractors in Abstract Parabolic Problems \publ Cambridge
University Press \publaddr Cambridge \yr 2000
\endref
\ref\no \dfa{F} \by E. Feireisl \paper Attractors for
semilinear damped wave equations on $\R^3$ \jour Nonlinear
Analysis \vol 23\yr 1994\pages 187--195
\endref
\ref\no
\dfa{F1} \bysame \paper Asymptotic behaviour and attractors
for semilinear damped wave equations with a supercritical
exponent \jour Proc. Roy. Soc. Edinburgh \vol 125A\yr
1995\pages 1051--1062
\endref
\ref\no \dfa{FY} \by D. Fall
and Y. You \paper Global attractors for the damped
nonlinear wave equation in unbounded domain \jour
Proceedings of the Fourth World Congress of Nonlinear Analysts\yr
2004\toappear
\endref
\ref\no \dfa{GT} \by J. M. Ghidaglia
and R. Temam \paper Attractors for damped nonlinear
hyperbolic equations \jour J. Math. Pures Appl. \vol 66\yr
1987\pages 273--319
\endref
\ref\no
\dfa{Go} \by J. A. Goldstein \book Semigroups of Linear Operators and applications \publ Oxford University Press
\publaddr New York \yr 1985
\endref
\ref\no \dfa{Ha} \by J. Hale
\book Asymptotic Behavior of Dissipative Systems \publ
American Mathematical Society \publaddr Providence \yr 1988
\endref
\ref\no \dfa{HR} \by J. Hale and G. Raugel \paper
Upper semicontinuity of the attractor for a singularly
perturbed hyperbolic equation \jour J. Differential
Equations  \vol 73\yr 1988\pages 197--214
\endref
\ref\no
\dfa{La1} \by O. Lady\v zenskaya \book The Boundary Value
Problems of Mathematical Physics \publ Springer-Verlag
\publaddr New York \yr 1985
\endref
\ref\no \dfa{La} \by O.
Lady\v zenskaya \book Attractors for Semigroups and
Evolution Equations \publ Cambridge University Press
\publaddr Cambridge \yr 1991
\endref
\ref\no \dfa{M} \by J.
E. Metcalfe \book Global Strichartz Estimates for Solutions
of the Wave Equation Exterior to a Convex Obstacle\bookinfo
PhD dissertation \publ Johns Hopkins University \publaddr
Baltimore \yr 2003
\endref
\ref\no \dfa{MRW} \by I. Moise, R. Rosa and X. Wang \paper
Attractors for non-compact semigroups via energy equations \jour Nonlinearity  \vol 11\yr 1998\pages 1369--1393
\endref
\ref\no \dfa{PR} \by M. Prizzi
and K. P. Rybakowski \paper Attractors for singularly
perturbed hyperbolic equations on unbounded domains
\paperinfo in preparation
\endref
\ref\no\dfa{Ra}
\by G. Raugel\paper Global attractors in partial differential equations  \inbook Handbook of dynamical systems, Vol. 2\pages  885--982\publ North-Holland\publaddr Amsterdam\yr 2002
\endref
\ref\no \dfa{SS} \by H.
F. Smith and C. D. Sogge \paper On Strichartz and
eigenfunction estimates for low regularity metrics \jour
Mathematical Research Letters\vol 1\yr 1994\pages 729--737
\endref
\ref\no \dfa{W} \by B. Wang \paper Attractors for
reaction-diffusion equations in unbounded domains \jour
Physica D\vol 179\yr 1999\pages 41--52
\endref
\endRefs
\enddocument